\documentclass[10pt,twoside]{article}

\usepackage{amsbsy,amsfonts,amsmath,amssymb,mathrsfs,diagrams,epsfig,eucal}

\def\itemm#1{\item[\indent {\rm (#1)}]}
\def\itemn#1{\item[\hspace{0.6mm} {\rm (#1)}]}
\def\qed{\hfill $\square$ \vspace{5mm}}
\def\titre#1{\par \vspace{2mm}\par \noindent {\bf {#1}} \par \vspace{2mm}\par}
\renewcommand{\ge}{\geqslant}
\renewcommand{\le}{\leqslant}
\def\twoquotient{/\!\!/}
\def\red{{\rm red}}
\def\zmod#1{\mathbb{Z}/#1\mathbb{Z}}
\def\et{{\rm \acute{e}t}}
\def\longto{\longrightarrow}
\def\isomto{\stackrel{\sim}{\longto}}

\newtheorem{definition}{Definition}[subsection]
\newenvironment{defi}{\begin{definition} \rm}{\end{definition}}

\newtheorem{prop}[definition]{Proposition}

\newtheorem{lemm}[definition]{Lemma}
\newtheorem{coro}[definition]{Corollary}
\newtheorem{theo}[definition]{Theorem}
\newtheorem{remark}[definition]{Remark}
\newenvironment{rema}{\begin{remark} \rm}{\end{remark}}
\newtheorem{remarks}[definition]{Remarks}
\newenvironment{remas}{\begin{remarks} \rm}{\end{remarks}}
\newtheorem{example}[definition]{Example}

\newtheorem{examples}[definition]{Examples}

\newtheorem{nothing}[definition]{$\!\!$}
\newenvironment{noth}{\begin{nothing} \rm}{\end{nothing}}
\newenvironment{proo}{{\flushleft \bf Proof :}}{\hfill $\square$ \vspace{5mm}}

\newtheorem{definition*}{Definition}[section]
\newenvironment{defi*}{\begin{definition*} \rm}{\end{definition*}}
\newtheorem{prop*}[definition*]{Proposition}
\newtheorem{lemm*}[definition*]{Lemma}
\newtheorem{coro*}[definition*]{Corollary}
\newtheorem{theo*}[definition*]{Theorem}
\newtheorem{remark*}[definition*]{Remark}
\newenvironment{rema*}{\begin{remark*} \rm}{\end{remark*}}
\newtheorem{examples*}[definition*]{Examples}
\newenvironment{exams*}{\begin{examples*} \rm}{\end{examples*}}

\DeclareMathOperator{\Hom}{Hom}
\DeclareMathOperator{\Isom}{Isom}
\DeclareMathOperator{\Aut}{Aut}
\DeclareMathOperator{\Id}{id}
\DeclareMathOperator{\Spec}{Spec}
\DeclareMathOperator{\Pic}{Pic}
\DeclareMathOperator{\Resul}{Res}
\DeclareMathOperator{\GL}{GL}
\DeclareMathOperator{\PGL}{PGL}
\DeclareMathOperator{\PSL}{PSL}
\DeclareMathOperator{\tr}{tr}

\def\clA{{\cal A}} \def\clL{{\cal L}} \def\clM{{\cal M}}
\def\clO{{\cal O}} \def\clP{{\cal P}} \def\clQ{{\cal Q}}

\def\bbA{{\mathbb A}} \def\bbD{{\mathbb D}} \def\bbF{{\mathbb F}}
\def\bbG{{\mathbb G}} \def\bbP{{\mathbb P}}
\def\bbQ{{\mathbb Q}} \def\bbZ{{\mathbb Z}}

\def\fre{\mathfrak{e}} \def\frf{\mathfrak{f}}
\def\frg{\mathfrak{g}} \def\frp{\mathfrak{p}}
\def\frA{\mathfrak{A}} \def\frS{\mathfrak{S}}

\newarrow{Dashto}{}{dash}{}{dash}{>}
\newarrow{Into}{C}{-}{}{-}{>}

\begin{document}

\begin{center}
{\bf \Large The stack of Potts curves and its fibre at a prime of wild}
\par \vspace{2mm}\par
{\bf \Large ramification}
\par \vspace{10mm}\par
{\bf Matthieu Romagny}
\par \vspace{2mm}\par
\emph{Stockholms Universitet, Matematiska institutionen,
106 91 Stockholm, Sweden}

\emph{e-mail address: romagny@matematik.su.se}
\end{center}

\par \vspace{12mm}\par

{\def\thefootnote{\relax}
\footnote{ \hspace{-6.8mm}
}
}

\par \vspace{1mm}\par
\footnotesize
\noindent {\bf Abstract}: In this note we study the modular properties of a
family of cyclic coverings of ${\bbP}^1$ of degree $N$, in all odd characteristics.
We compute the moduli space of the corresponding algebraic stack over $\bbZ[1/2]$,
as well as the Picard groups over algebraically closed fields. We put special
emphasis on the study of the fibre of the stack at a prime of wild ramification;
in particular we show that the moduli space has good reduction at such a prime.
\normalsize

\par \vspace{10mm}\par

We are interested here in algebraic curves called \emph{$N$-state Potts curves},
previously studied here and there in the literature but not from a modular and
arithmetic viewpoint. They provide an interesting example (quite unique in fact)
of a stack of Galois covers where the study can be pursued quite far even at
the primes of wild ramification. We discover interesting phenomena, as well as
explanations and commentaries concerning the deformation theory of wild covers
of curves (see~\cite{BM1}).

In the modular theory of covers of curves, much is known in the tame
case where the characteristic $p\ge 0$ of the base field does not divide the
indices of ramification: these covers form a  smooth algebraic stack, and it
is known how to compactify it with stable covers. By contrast, for example
for a Galois cover of group $G$ whose order is divisible by $p$, the
(uni)versal deformation ring is very nasty (in general not even reduced),
meaning that the corresponding classifying stack is far from being smooth.
Moreover, there exist smooth $G$-covers in characteristic~0 that do not
have good (i.e. smooth) reduction to characteristic~$p$, and vice versa
there exist smooth $G$-covers in char.~$p$ that do not lift to char.~0.
Giving a bridge between these characteristics usually means studying
objects over a valuation ring of mixed characteristics $(0,p)$, and how they
specialize or generize. In the example of Potts curves everything can even
be done over $\bbZ[1/2]$, and we will be able to decribe quite precisely the
behaviour at a prime of wild ramification.

\par \vspace{2mm}\par

Let us now describe in more detail the results of the article. Given an odd
integer $N\ge 3$, an $N$-state Potts curve (or $N$-Potts curve) is by
definition a smooth hyperelliptic curve of genus $N-1$, which is a cyclic
covering of $\bbP^1$ of degree $N$. We will simplify the treatment of
hyperellipticity by avoiding the prime $p=2$, and we will rather focus on
wild ramification at primes dividing $N$. Hence in all the article the schemes
and stacks considered are over $\bbZ[1/2]$. Thus the fibered category with
objects $N$-Potts curves is an algebraic stack $\clP_N$ over $\bbZ[1/2]$.
Among the main results of the article is the following computation
(theorems~\ref{espace_modulaire} and~\ref{espace_modulaire_anychar}):
\begin{itemize}
\item If $N$ is not prime, the coarse moduli space of $\clP_N$ is the
scheme $\bbA^1_*\otimes\bbZ[\frac{\zeta_{_N}+\zeta_{_N}^{-1}}{2},\frac{1}{2N}]$.
\item If $N=p$ is prime, the coarse moduli space of $\clP_p$ is the
scheme $\bbA^1_*\otimes\bbZ[\frac{\zeta_{_p}+\zeta_{_p}^{-1}}{2},\frac{1}{2}]$.
\end{itemize}
In this statement, $\bbA^1_*=\bbA^1-\{0\}$ is the punctured line and the
arithmetic rings that appear are subextensions of degree 2
of the ring of cyclotomic integers (see~\ref{remaMatrix}). Note that when $N$
is a composite integer, there does not exist Potts curves in characteristics
$p|N$. We see that, although the stacks (and moduli spaces) are connected for
arithmetic reasons, at (geometric) primes of tame ramification they split as
a sum of $\varphi(N)/2$ connected irreducible components, where $\varphi$
is the Euler function. On the contrary, when $N=p$ is prime,
we show that the moduli space "has good reduction" at $p$, that is to
say its formation commutes with the base change $\bbZ[1/2]\to \bbF_p$~:
\begin{itemize}
\item The coarse moduli space of $\clP_p\otimes \bbF_p$ is
$\bbA^1_*\otimes \frac{\bbF_p[z]}{z^{(p-1)/2}}$, the fibre
of the moduli space of $\clP_p$ at $p$.
\end{itemize}
The result (theorem~\ref{espace_modulaire_carp}) is even more precise
and shows that the map from $\clP_p\otimes \bbF_p$ to its moduli space
is \'etale, and so this stack is connected but
non-reduced with multiplicty $\varphi(p)/2$; its reduced part is smooth.
Hence we can interpret the non-reducedness as coming from the collision
of $\varphi(p)/2$ smooth components in characteristic 0 when we reduce
to characteristic $p$.

We also compute the Picard groups of the geometric fibres of $\clP_N\otimes k$
($k$ an algebraically closed field of characteristic $p\ne 2$).
The first fact to be noted is that these groups are finite at primes of
tame ramification, but not anymore at primes of wild ramification.
The second interesting point is that the nilpotents contribute a lot
to the Picard group of $\clP_p\otimes \overline{\bbF}_p$, which would
certainly not be the case for its moduli space $P$ because it is an
affine scheme, hence $\Pic(P_\red)=\Pic(P)$. The reason for this
is of course the presence of automorphisms. Here is the result
(theorems~\ref{PicardTame} and~\ref{PicardWild})~:
\begin{itemize}
\item If $(N,p)=1$ then the Picard group of any connected component
of $\clP_N\otimes k$ is isomorphic to $(\zmod{2})\times (\zmod{2N})$.
\item If $N=p$ then the Picard group of $\clP_p\otimes k$ is isomorphic
to $\zmod{2}\times (1+zA)$, where $A=\frac{k[z]}{z^{(p-1)/2}}[X,\frac{1}{X}]$
is the affine ring of the moduli space, and $1+zA\subset A^\times$ is
a multiplicative subgroup.
\end{itemize}
Finally we mention that in the tame case, it is likely that a little
more work would give the expected results concerning the stack of
\emph{stable} $N$-Potts curves, namely, that its moduli space is
$\bbP^1\otimes\bbZ[\frac{\zeta+\zeta^{-1}}{2},\frac{1}{2N}]$. We will
say no more than a word about this in~\ref{CCS}. Concerning the stack
of stable $p$-Potts curves over $\bbZ[1/2]$, similar questions arise
but here the problems are of course more complicated. It is clear that
the work done in the present article makes it tempting to ask what
would stable $p$-Potts curves look like; if there is a 1-dimensional stack
of stable $N$-Potts curves even over characteristics $p|N$; what are the
stable reductions of $N$-Potts curves in characteristic $p$, when $p$
divides a non-prime $N$... All this is left aside for the time being.

\par \vspace{2mm}\par

Here is a short overview of the organization of the article. In the text,
the order of apparition of the results is actually rather different from
the one above. The first section contains preliminaries on
finite subgroups of the projective linear group ${\rm PGL}_2$. The second
section contains the computation of the automorphism groups of Potts curves.
Here we make the essential observation that the good understanding of the
stack $\clP_N$ over $\bbZ[1/2N]$ is via the classical Hurwitz description
of branched covers by the shape of the ramification. This leads us to start
from a different definition for Potts curves, and of course we eventually
show that the two coincide. In the third section, we treat the case of a
composite~$N$: we compute the moduli space of $\clP_N$ and the Picard group
of its fibres. In the fourth section we extend the construction of the moduli
space to the case of the stack $\clP_p$ with $p$ prime. At last the fifth
section is devoted to the study of the fibre of $\clP_p$ at $p$.

\par \vspace{2mm}\par

\emph{Conventions}. We will consider that every positive integer is prime to~0, so
as not to make repeated particular cases when speaking of primality of an
integer with the characteristic of a field. The Euler function is denoted
$\varphi$ as usual. Finally, in a module $M$ over a commutative ring $A$,
we will denote by $m\propto m'$ the equality up to multiplication by an
invertible element of $A$.

\emph{Acknolewdgements}. It has been a long awaited moment the time to thank here
my advisor Jos\'e Bertin, who followed from the beginning the story of Potts curves,
with a lot of energy and disponibility. Also I express warm thanks to Laurent
Moret-Bailly and Ariane M\'ezard, who read former versions of this work. Their
thorough reading and numerous comments were the source of many improvements.
Finally I acknowledge a grant from the EAGER network, hospitality of Stockholm
University and the Institut Fourier in Grenoble where most of this work was done.

\section{Preliminaries on Dickson's theorem}

There is no pretence to originality in the contents of this section, but the
results stated here could not be found in the literature. Because of their
simplicity, proofs are sometimes elliptical, if not omitted. The reader may
without prejudice skip this section and go straight to \S~\ref{section_CP_autom},
refering to the results below when necessary.

In section~\ref{section_CP_autom}, in the course of the computation of the
automorphism group of Potts curves over an algebraically closed field $k$ of
characteristic $p\ne 2$ (see~\ref{theo_autom} and~\ref{automorphismes_car_p}),
we will have to handle certain finite subgroups of $\PGL_2(k)$. Recall
from~\cite{Su}, chap.~3, th.~6.17 that Dickson's theorem for $\PGL_2$ takes
the following shape:

\begin{theo*}[Dickson] \label{Dickson}
Let $k$ be an algebraically closed field of characteristic $p\ne 2$. Any finite
subgroup $G\subset \PGL_2(k)$ is isomorphic to a subgroup among the following list: \\
$\bullet$ If $p$ is prime to $|G|$,
\begin{trivlist}
\itemm{1} cyclic group, dihedral group, symmetric $\frS_4$,
alternating $\frA_4$ or $\frA_5$,
\end{trivlist}
$\bullet$ If $p$ divides $|G|$,
\begin{trivlist}
\itemm{2} $G=Q\rtimes C$ a semi-direct product of a normal, elementary abelian,
$p$-Sylow $Q$ by a cyclic group of order prime to $p$,
\itemm{3} $\frA_5$ if $p=3$,
\itemm{4} $\PSL_2(\bbF_q)$ or $\PGL_2(\bbF_q)$ for $q=p^s$, $s\ge 1$ integer. \qed
\end{trivlist}
\end{theo*}

In this section we classify \emph{conjugation} classes instead of merely
\emph{isomorphism} classes. We define an equivalence relation in $k^\times$ by
$x'\sim x$ $\Leftrightarrow$ $x'\in\{x,x^{-1}\}$. The
corresponding class is denoted $[x]$, and the mapping $[x]\mapsto
x+x^{-1}$ is a set-theoretic injection $k^\times/\!\sim\,\hookrightarrow k$.
Finally let $\mu_n^*\subset k^\times$ be the set of primitive $n$-th roots of
unity (e.g. $\mu_p^*=\{1\}$).

\begin{prop*} \label{automorphismes_PGL2}
Let $A\in \PGL_2(k)$ be an automorphism of finite order $n>1$. Then,
\begin{trivlist}
\itemn{i} Either $n$ is prime to $p$, or equal to $p$.
\itemn{ii} As an automorphism of $\bbP^1_k$, $A$ is conjugated to $x\mapsto \zeta x$,
for some $\zeta\in\mu_n^*$, when $(n,p)=1$, and to $x\mapsto x+1$, when $n=p$.
\itemn{iii}
The set $\mu_n^*/\!\sim$ classifies conjugation classes of elements of order $n$,
and more precisely,
\end{trivlist}
$$
{\rm ord}(A)=n \; \Leftrightarrow \; \mbox{ there exists }
[\zeta]\in\mu_n^*/\!\sim \mbox{ such that }
(\zeta+\zeta^{-1}+2)\det(A)-\tr(A)^2=0
$$
\end{prop*}

\begin{proo}
We work with a representative in $\GL_2(k)$, whose eigenvalues are given
by the characteristic polynomial as $\lambda^{\pm}=\frac{1}{2}(\tr(A)\pm\delta)$
with $\delta^2=\tr(A)^2-4\det(A)$. In $\PGL_2(k)$, only the class
(for the relation~$\sim$) of the ratio $\zeta:=\lambda^+/\lambda^-$
is well-determined. 

We have $[\zeta]=1$ if and only if, up to homothety, $A$ is conjugated to
a unipotent matrix, i.e. $n=p$ and, as a homography, $A$ is conjugated
to $x\mapsto x+1$. We have $[\zeta]\ne 1$ if and only if $A$ is conjugated
to the diagonal matrix diag($\lambda^+,\lambda^-$), and then $A$ has order
$n$ if and only if $(n,p)=1$ and $[\zeta]\in \mu_n^*/\!\sim$. As a homography,
$A$ is conjugated to $x\mapsto \zeta x$.

For the claim in (iii), one needs just checking that
$$\zeta+\zeta^{-1}=\frac{\tr(A)+\delta}{\tr(A)-\delta}+
\frac{\tr(A)-\delta}{\tr(A)+\delta}=\frac{\tr(A)^2}{\det(A)}-2$$
\end{proo}

\begin{exams*} \label{examsOrders}
As particular cases of~\ref{automorphismes_PGL2}(iii), an element $A\in \PGL_2(k)$
has order~2 (resp. order 3, 4, 6 or $p$) if and only if $\tr(A)^2=i\det(A)$ for $i=0$
(resp. $i=1$, 2, 3 or 4).
\end{exams*}

\begin{rema*} \label{remaMatrix}
Let $\zeta$ be a primitive $n$-th root of unity (say as a complex
number) and $\Phi_n$ the $n$-th cyclotomic polynomial (of degree $\varphi(n)$).
For future use, we observe that $\frac{1}{2}(\zeta+\zeta^{-1})$ is integral
over $\bbZ[\frac{1}{2}]$. Indeed its minimal polynomial over $\bbQ$ is
$2^{-\varphi(n)/2}\psi_n(2t)$ where $\psi_n\in\bbZ[t]$ is the monic polynomial
such that $\Phi_n(t)=t^{\varphi(n)/2}\psi_n(t+t^{-1})$. According to (iii) of
the proposition we define an automorphism of the projective line over the
spectrum of $\bbZ[\frac{\zeta+\zeta^{-1}}{2},\frac{1}{2n}]$, of exact order
$n$ on all the fibres, by the following matrix:
$$
M_\zeta=\left(
\begin{array}{cc}
\zeta+\zeta^{-1} & \zeta+\zeta^{-1}-2 \\
1 & 2
\end{array}
\right)
$$
\end{rema*}

\par \vspace{2mm}\par

\begin{coro*} \label{orders_in_PGL2(Fq)}
Let $q=p^s$ for some $s\ge 1$. Then the order of an element
$A\in\PGL_2(\bbF_q)$ divides $q-1$, $q+1$, or $p$.
\end{coro*}

\begin{proo}
Take $A$ an element of order $n$ prime to $p$. By proposition~\ref{automorphismes_PGL2},
there exists $\zeta\in\mu_n^*$ algebraic over $\bbF_q$, of degree at most 2,
such that $A$ is conjugated to $x\mapsto \zeta x$. If
$\zeta\in\bbF_q$, we have $\zeta^{q-1}=1$. Else, the minimal polynomial of $\zeta$
over $\bbF_q$ is $P=X^2+(2-\frac{\tr(A)^2}{\det(A)})X+1$. But $P$ can
also be written $P=X^2-(\zeta+\zeta^q)X+\zeta^{q+1}$ with the Frobenius ${\rm Fr}(x)=x^q$,
generating ${\rm Gal}(\bbF_{q^2}/\bbF_q)=\zmod{2}$. Hence $\zeta^{q+1}=1$ and we are done.
\end{proo}

\begin{coro*} \label{fixed_points}
There are $q^2-1$ elements of order $p$ in ${\rm PGL}_2(\bbF_q)$, and they all
belong to ${\rm PSL}_2(\bbF_q)$. Moreover the set of fixed points of all order $p$
elements of ${\rm PGL}_2(\bbF_q)$, acting on $\bbP^1(k)$, is $\bbP^1(\bbF_q)$.
\end{coro*}

\begin{proo}
Simple calculations.
\end{proo}



In the sequel we use the concise notation $<\alpha(x)>$ for the
subgroup generated by a homography $\alpha\in{\rm PGL}_2(k)$.

\begin{coro*} \label{cyclic_dihedral}
The cyclic and dihedral subgroups of $\PGL_2(k)$ are conjugated to:
\begin{trivlist}
\itemn{i}
$\left\{
\begin{array}{ll}
{\sf C}_n=\,<\zeta\,x> & \mbox{for } (n,p)=1 \;(\mbox{any }\zeta\in\mu_n^*) \\
{\sf C}_p=\,<x+1>     & \mbox{for } n=p.
\end{array}
\right.$
\itemn{ii}
$\left\{
\begin{array}{ll}
{\sf D}_n=\,<\zeta\,x,\,\frac{1}{x}> & \mbox{for } (n,p)=1 \;(\mbox{any }\zeta\in\mu_n^*) \\
{\sf D}_p=\,<x+1,-x>                & \mbox{for } n=p.
\end{array}
\right.$ \qed
\end{trivlist}
\end{coro*}

\begin{coro*} \label{S4}
The subgroups of $\PGL_2(k)$ isomorphic to $\frS_4$ ($p\ne 2,3$) are conjugated to
$$ {\sf S}_4\:=\:<ix,\frac{x+1}{x-1}>\:\simeq \: <(1234),\,(12)>\;.$$
\end{coro*}

\begin{proo}
Let $\nu:\frS_4 \stackrel{\sim}{\longrightarrow} G$ be
an isomorphism with values in a subgroup of ${\rm PGL}_2(k)$;
$a=\nu (1234)$ and $b=\nu (12)$ generate $G$. Up to conjugation,
$a(x)=ix$. A priori the involution $b$ can be written $b(x)=\frac{rx+s}{tx-r}$;
using the fact that $ab=\nu (134)$ has order 3, we get $r^2=st$
(see~\ref{examsOrders}). Conjugation by $\phi(x)=rx/t$ leaves $a$
invariant while $b$ maps to the desired $x\mapsto \frac{x+1}{x-1}$.
To complete the proposition, one checks that these two elements
generate a subgroup isomorphic to $\frS_4$.
\end{proo}

\begin{coro*} \label{PGL2/PSL2}
The subgroups isomorphic to ${\rm PGL}_2(\bbF_q)$ are all conjugated to the
"standard" ${\rm PGL}_2(\bbF_q)$ corresponding to the field inclusion 
$\bbF_q\hookrightarrow k$; the same result holds for ${\rm PSL}_2(\bbF_q)$.
\end{coro*}

\begin{proo}
The standard $\PGL_2(\bbF_q)$ is generated by the following three elements:
$$
e(x)=x+1,\quad f(x)=ux,\quad g(x)=1/x
$$
where $u$ is any generator of the multiplicative group $\bbF_q^*$. Setting
$m=\frac{q-1}{p-1}$ then $v=u^m$ is a generator of $\bbF_p^*$; in particular
$v$ is an integer modulo $p$. Then we have a relation $e^vf^m=f^me$, it is
just the homography $x\mapsto vx+v$. Also it is immediate that $n:={\rm ord}(eg)$
is prime to $p$.

Now let $G$ be a subgroup of $\PGL_2(k)$, and
$\nu:\PGL_2(\bbF_q)\to G$ be an isomorphism. Denote $\fre=\nu(e),\frf=\nu(f),\frg=\nu(g)$
so that $G=<\fre,\frf,\frg>$. As $\frf,\frg$ generate a dihedral group,
by the above result~\ref{cyclic_dihedral}, with a first conjugation we can
suppose that $\frf=f$ and $\frg=g$. The above relations in $\PGL_2(\bbF_q)$
yield:
$$
(\dag) \:  \fre^p=1 \qquad
(\dag\!\dag) \:   \fre^v\frf^m=\frf^m\fre \qquad
(\dag\!\dag\!\dag)  \: (\fre\frg)^n=1
$$
Let us write $\fre(x)=\frac{ax+b}{cx+d}$, then $(\dag)$ reads $(a+d)^2=4(ad-bc)$
by~\ref{automorphismes_PGL2}. By induction,
$$
\fre^k(x)=\frac{\left(\frac{k+1}{2}a-\frac{k-1}{2}d\right)x+kb}
{kcx-\frac{k-1}{2}a+\frac{k+1}{2}d}
$$
We then write $(\dag\!\dag)$ explicitly and obtain $a=d$, and $c=0$. At this point,
$\fre(x)=x+b/a$. Then by $(\dag\!\dag\!\dag)$ and~\ref{automorphismes_PGL2} there
exists $\lambda\in \mu_n^*\subset \bbF_{q^2}^*$
such that $(1+\lambda)^2 a^2=-\lambda b^2$. As $n$ divides $q-1$ or $q+1$
by~\ref{orders_in_PGL2(Fq)}, we have $\lambda+\lambda^{-1}\in\bbF_q$ and so
$\beta:=\left(\frac{b}{a}\right)^2=-(\lambda+\lambda^{-1}+2)\in\bbF_q$.
Finally we apply a conjugation by $\phi:x\mapsto \frac{b}{a}x$, then $G=<~x+b/a,ux,1/x~>$
is mapped to $<x+1,ux,\frac{1}{\beta x}>={\rm PGL}_2(\bbF_q)$ as announced.
The case of ${\rm PSL}_2(\bbF_q)$ is similar.
\end{proo}

\section{Potts curves and their automorphisms} \label{section_CP_autom}

In this section we describe Potts curves defined over an algebraically
closed field, before looking at families (hence moduli) in the rest of
the article. Let us make precise definitions before we start: in all what
follows, by a curve over a scheme $S$ we will mean a proper, flat morphism
$f:C\to S$ whose fibres are projective, geometrically connected and
one-dimensional; also, except in~\ref{CCS}, they will be assumed to be smooth.
We recall that:

\begin{defi*} \label{def0}
Let $N\ge 3$ be an odd integer. An $N$-Potts curve is a smooth hyperelliptic
curve of genus $N-1$, which is a cyclic covering of $\bbP^1$ of degree $N$.
\end{defi*}

In subsection~\ref{CMH} where $N$ is prime to $p$,
we have three main goals: showing the equivalence between two different
definitions of Potts curves~(\ref{def0<=>def2.1.1}), giving a modular
invariant~(\ref{j-isom}), and computing automorphism groups~(\ref{theo_autom}
and \ref{coroAutMod}). In subsection~\ref{CS} where $p$ divides $N$,
we also define an invariant and compute the automorphis groups (it turns
out to be simpler).

\subsection{Tame case} \label{CMH}

In this case, the curves we are dealing with can be described like in the classical
Hurwitz setting, i.e. as maps to $\bbP^1$ with prescribed ramification.
This is our starting point; it leads us to change definition~\ref{def0}.

Recall that if a finite group $G$ acts faithfully on a smooth curve $C$
over a field of characteristic~$p$, such that $|G|$ and $p$ are coprime,
then the stabilizer of a fixed point is cyclic and its natural
representation in the cotangent space of the point is faithful.
This gives rise to the \emph{Hurwitz ramification datum}, i.e. the
list of all the corresponding characters at the fixed points.

Let $N\ge 3$ be an odd integer. Set $G=\zmod{N}$, and assume that a
generator for $G$, denoted~$g$, has been chosen once for all. In the
character group $\widehat{G}\simeq G$, we define an equivalence relation
by $\chi'\sim\chi$ if and only if $\chi'\in\{\chi,\chi^{-1}\}$. Denote
by $[\chi]$ the corresponding class.

\begin{defi} \label{def_Potts_curve}
Let $k$ be an algebraically closed field of characteristic $p$ prime to $2N$.
\begin{trivlist}
\itemn{i} An $N$-Potts curve of type $[\chi]$ over $k$ is a curve $C$
together with a faithful action $\rho:G\hookrightarrow \Aut_k(C)$
such that $C/G$ has genus 0, with four ramification points all with
stabilizer equal to $G$, and Hurwitz ramification datum
$\{\chi,\chi,\chi^{-1},\chi^{-1}\}$. An $N$-Potts curve is an $N$-Potts
curve of type $[\chi]$ for some~$[\chi]$.
\itemn{ii} An isomorphism between two Potts curves $C,C'$ is a $G$-equivariant
isomorphism of algebraic curves $\varphi:C\stackrel{\sim}{\longrightarrow} C'$.
\end{trivlist}
\end{defi}

\begin{remas} \label{csq_def} \rm
\begin{trivlist}
\itemn{i} The action $\rho$ being determined by $\sigma=\rho(g)$,
a Potts curve will be denoted $(C,\sigma)$. If $(C,\sigma)$ and $(C',\sigma')$
are isomorphic $N$-Potts curves then $[\chi]=[\chi']$.
\itemn{ii} Of course, from now up to~\ref{def0<=>def2.1.1} it is this
definition that applies rather than definition~\ref{def0}.
\itemn{iii} If $\sigma$ acts via $\chi$ on the cotangent
space $m_a/m_a^2$, then $\chi$ is determined by the root of
unity $\zeta:=\chi(\sigma)$. Hence, the quantity $\zeta+\zeta^{-1}$
is attached to $C$, and equivalent to~$[\chi]$.
\end{trivlist}
\end{remas}

\titre{Hyperellipticity}

Let $(C,\sigma)$ be an $N$-Potts curve. By the Riemann-Hurwitz formula,
we get the genus $g(C)=N-1$. Notice that a birational equation can
easily be drawn from the definition: by Kummer theory, the function field
$k(C)$ is generated, as an extension of the rational function field
$k(x)$, by a single $t\in k(C)$ such that $t^N\in k(x)$, i.e.
\begin{equation} \label{EQ0}
t^N=\frac{(x-a)(x-b)}{(x-c)(x-d)}
\end{equation}
according to the ramification. Substituting $t/(x-c)(x-d)$ to $t$ we get
$$
t^N=(x-a)(x-b)(x-c)^{N-1}(x-d)^{N-1}
$$
where $\sigma$ acts by $t\mapsto \zeta t$ for some $\zeta$. Now in
${\rm Aut}_k(\mathbb{P}^1)$ there is one and only one subgroup isomorphic
to $\zmod{2} \times \zmod{2}$, generated by involutions that interchange
the four points $a,b,c,d$ (see~\cite{R}, \S~1, Lemma~2), namely
$$
\begin{array}{ccc}

\tau_0:\;\left\{ 
\begin{array}{l}
a \leftrightarrow b \\
c \leftrightarrow d
\end{array}
\right.
&
\qquad\mu_0:\;\left\{ 
\begin{array}{l}
a \leftrightarrow c \\
b \leftrightarrow d
\end{array}
\right.
&
\qquad\tau_0\mu_0=\mu_0\tau_0:\;\left\{ 
\begin{array}{l}
a \leftrightarrow d \\
b \leftrightarrow c
\end{array}
\right.

\end{array}
$$


In order to make the link between~\ref{def_Potts_curve} and definition~\ref{def0}
we must build a hyperelliptic involution from $\tau_0$ and for that
we need to recall the following classical construction:

\par \vspace{4mm}\par

\begin{noth} \label{X/G_<=>_(Y,L,s)} \rm
It is known how to describe a tame cyclic covering of curves (or even,
families of curves) in terms of invertible sheaves on the base. A perfect
exposition is recalled in \cite{ArV}, so we just give a sketch here.
For $n\ge 2$, let $S$ be a scheme with $n\in {\cal O}_S^\times$, and
$\zeta\in {\cal O}_S$ a primitive $n$-th root of unity. Let $X\to S$
be a smooth curve, and $\sigma$ an $S$-automorphism of order~$n$.
By smoothness the quotient morphism $f:X\to Y=X/\sigma$ is finite flat.
By the assumption of invertibility of $n$, there is a decomposition
$f_*{\cal O}_X=\oplus_{j=0}^{n-1}\,\clL_j$ with $\clL_j$ equal to
the $\zeta^j$-eigenspace for the action of $\sigma$.
The $\clL_j$ are invertible sheaves, and the multiplication in $f_*{\cal O}_X$
gives injective maps $\clL_i\otimes\clL_j\to\clL_{i+j}$ ($i+j$ is read
modulo $n$). In particular,
as $\sigma^n={\rm id}$, we get $\clL_{_1}^n\simeq {\cal O}_Y(-D)$
where $D$ is the effective Cartier branch divisor of $f$.

Conversely, given $\clL=\clL_{_1}$ and a global section $s$ whose divisor
of zeroes is $D$ (so $s$ is determined up to a global invertible section),
we can reconstruct the $\clL_j$ and endow $\clA=\oplus_{j=0}^{n-1}\,\clL_j$
with a product mapping $(\ell,\ell')\in \clL_j\times \clL_{n-j}$ to
$s\ell\ell'\in {\cal O}_Y$, and we recover $X=\Spec(\clA)$.

As a conclusion, we can consider the datum of an $S$-curve $X$ with an
automorphism of order $n$ as being equivalent, up to isomorphism, to
that of a triple $(Y,\clL,s)$ where $\clL\in {\rm Pic}(Y)$ and $s$ is a
global section of $\clL^{-n}$. Furthermore there is an obvious fonctoriality
in $(Y,\clL,s)$: for a map $(Y',\clL',s')\to (Y,\clL,s)$ given by an affine
$S$-morphism $\alpha:Y'\to Y$ and a map $\beta:\clL\to \alpha_*\clL'$
respecting the sections, there is an induced morphism $h:X'\to X$ with
$\sigma h=h\sigma'$. \qed
\end{noth}

\par \vspace{2mm}\par

Back to our situation, the action of $\sigma$ on the Potts curve $C$
is described by $\clL={\cal O}(-2)$ and $s=(X-a)(X-b)(X-c)^{N-1}(X-d)^{N-1}$.
Clearly $\tau_0$ gives an automorphism of the triple $(\bbP^1,\clL,s)$,
hence lifts to an automorphism $\tau:C\to C$.

Moreover, $\tau$ satisfies $\tau\sigma\tau^{-1}=\sigma$. Also we have that
$\tau^2=\sigma^j$ for some $j\in \zmod{N}$, because it induces the identity
on $\bbP^1$. But 2 being invertible modulo $N$, we may write $j=2k$, and then,
changing $\tau$ into $\tau\sigma^{-k}$ if necessary, we can assume that $\tau^2=1$.

As for $\mu_0$, things are slightly different because it exchanges $a,c$ and
$b,d$. Let $s'$ be the section $(X-c)(X-d)(X-a)^{N-1}(X-b)^{N-1}$ of
$\clL={\cal O}(-2)$, then $\mu_0$ gives a map between $(\bbP^1,\clL,s)$
and $(\bbP^1,\clL,s')$. The last triple gives rise to the same curve $C$
of course, but with the automorphism $\sigma^{-1}$. So $\mu_0$ lifts to
$\mu:C\to C$ such that $\sigma^{-1}\mu=\mu\sigma$. As above, we may change $\mu$
so as to have $\mu^2=1$; then $\mu$ and $\sigma$ generate a group isomorphic
to the dihedral group $\bbD_N$.

\par \vspace{2mm}\par

\begin{prop} \label{Potts=>hyperelliptique}
$C$ is hyperelliptic, and $\tau$ is the hyperelliptic involution.
\end{prop}

\begin{proo}
Using the fact that $N$ and 2 are coprime, it is immediate that a
point in $C$ is fixed by $\tau$ if and only if its image in $C/\sigma$
is fixed by $\tau_0$. But the supports of the ramification loci for
the quotients by $\tau$ and $\sigma$ (i.e. their fixed points) are
disjoint, because their images in $C/\sigma$ are already.
Hence we get $2N$ fixed points for $\tau$, namely all the preimages
of the fixed points of $\tau_0$. They form two $\sigma$-orbits.
Applying the Riemann-Hurwitz formula to the quotient $C\to C/\tau$ of degree 2:
$$
2(N-1)-2=2\,(2g_{C/\tau}-2)+2N
$$
we have $g_{C/\tau}=0$, as desired.
\end{proo}

\begin{prop} \label{def0<=>def2.1.1}
Definitions~\ref{def0} and \ref{def_Potts_curve} are equivalent.
\end{prop}

\begin{proo}
It only remains to prove the implication "\ref{def0} $\Rightarrow$ \ref{def_Potts_curve}".
But once again, using oddness of $N$, a point in $C$ is fixed by $\sigma$ if
and only if its image in $C/\tau$ is fixed by $\sigma_0$ (the morphism
induced from $\sigma$), and then the stabilizers are equal. An automorphism
of $\bbP^1$ of order $N$ has two fixed points with full stabilizer,
and ramification characters inverses to each other. As the ramification loci
for $\sigma$ and $\tau$ are disjoint, the two points lifted in $C$ give four
fixed points with stabilizer $G=\zmod{N}$ as in~\ref{def_Potts_curve}, and
with the expected Hurwitz ramification datum.
\end{proo}

\begin{rema} \label{chgt_definition} \rm
By the way, in \cite{R}, the definition chosen for Potts curves is a mix
between ours: there, an $N$-Potts curve is a hyperelliptic curve of genus
$N-1$, with an order $N$ automorphism having exactly 4 fixed points.
\end{rema}

\titre{Automorphisms}

Using the quotient by $\tau$ we can get another affine equation for a Potts
curve $C$. As a matter of fact, $\sigma$ induces on $C/\tau\simeq \bbP^1$ an
automorphism conjugated to $x\mapsto \zeta x$. The branch locus of $\tau$
is composed of two orbits of $\sigma$, i.e.
$\{\zeta^j \alpha\}\cup\{\zeta^j \beta\}$ for $0\leq j\leq N-1$, for certain
$\alpha,\beta$ with $\alpha,\beta$ and 0 all distinct. The corresponding equation is
\begin{equation} \label{EQ2}
y^2=(x^N-\alpha^N)(x^N-\beta^N)=x^{2N}+Ax^N+B
\end{equation}
We can recover (\ref{EQ0}) with the choice of a rational parameter for the
conic $y^2=u^2+Au+B$ (e.g. with the coordinates $z=\frac{y+\sqrt{B}}{x^N}$
and $t=\frac{x}{(2\sqrt{B})^{1/N}}$, the equation is
$t^N=\frac{z-\lambda}{z^2-1}$ with $\lambda=\frac{-A}{2\sqrt{B}}$).
The automorphisms $\sigma,\tau,\mu$ have the following
expressions on model~(\ref{EQ2})~:
$$
\begin{array}{ccc}
\sigma(x,y)=(\zeta x,y) \quad; &
\tau(x,y)=(x,-y) \quad; &
\mu(x,y)=(\frac{B^{1/N}}{x},\frac{\sqrt{B}y}{x^N}) \;.
\end{array}
$$

\noindent Let us define a modular invariant $j=\frac{B}{A^2-4B}\ne 0$
for a curve with equation~(\ref{EQ2}).
As is expected,

\begin{prop} \label{j-isom}
Two Potts curves $(C,\sigma)$ and $(C',\sigma')$ of invariants $j$ and $j'$
are isomorphic if and only if $j=j'$ and $[\chi]=[\chi']$.
\end{prop}

\begin{proo}
Let $\varphi:C\to C'$ be an isomorphism with $\sigma'\varphi=\varphi\sigma$.
It induces a map $\tilde{\varphi}:C/\tau\to C'/\tau'$ on the quotients.
Moreover, proposition~\ref{automorphismes_PGL2} says that for an automorphism
of $\bbP^1$ of order prime to $p$, there is a unique coordinate $x$ on $\bbP^1$
(up to $\sim$) such that the automorphism is a homothety. So if we consider
equations of type~(\ref{EQ2}) for $C$ and $C'$, then $G$-equivariance reads
either $\tilde{\varphi}(\zeta x)=\zeta\tilde{\varphi}(x)$,
or $\tilde{\varphi}(\zeta x)=\zeta^{-1}\tilde{\varphi}(x)$.
In the first case we find $\tilde{\varphi}(x)=\lambda x$ for some $\lambda$,
whence $A'=\lambda^N A$, $B'=\lambda^{2N} B$ and $j=j'$.
In the second case we find $\tilde{\varphi}(x)=\lambda /x$ for some $\lambda$,
whence $A'=\lambda^N A/B$, $B'=\lambda^{2N}/B$ and $j=j'$.

Conversely, assume that $j=j'$ and $[\chi]=[\chi']$ with $C$ and $C'$ given
by an equation~(\ref{EQ2}). Then $\chi'=$ either $\chi$ or $\chi^{-1}$. If
$\chi'=\chi$ choose $\lambda$ such that $A'=\lambda^N A$ and $B'=\lambda^{2N} B$;
then $\varphi:(x,y)\mapsto (\lambda x,\lambda^N y)$ is a $G$-isomorphism from
$C$ to $C'$. If $\chi'=\chi^{-1}$ choose $\lambda$ such that $A'=\lambda^N A/B$
and $B'=\lambda^{2N}/B$; then
$\varphi:(x,y)\mapsto(\lambda /x,\sqrt{B}(\lambda /x)^Ny)$ answers the question.
\end{proo}

\begin{rema} \label{remaA1*}
Via $j$, there is a 1-1 correspondance between isomorphism classes of
$N$-Potts curves and the sum of $\varphi(N)/2$ copies of $\bbA^1-\{0\}$.
Indeed, for fixed $[\chi]$ and $j\ne 0$, a curve with invariant $j$ is
given by the equation $y^2=x^{2N}+(1+4j)x^N+j(1+4j)$ if $j\ne -1/4$, or
by $y^2=x^{2N}-1$ if $j= -1/4$.
\end{rema}

We are now in position to compute the automorphism groups of Potts curves,
using the results of section 1. On the way, we provide a correction to
\cite{R}, section~2, prop.~5 where the case of $j=-1/4$ was forgotten.

\begin{theo} \label{theo_autom}
Let $N\ge 3$ be an odd integer, and $k$ an algebraically closed field
of characteristic $p$ prime to $2N$. Let $(C,\sigma)$ be an $N$-Potts curve.
Then the automorphism group of the curve $C$ (alone) is the following:
\begin{trivlist}
\itemn{i} If $p\ne 3,5$, $N=3$, $j=-1/54$, then
${\rm Aut}_{k}(C)=\widetilde{\frS_4}$,
the representation group of $\frS_4$ where the elements corresponding to
transpositions have order 2, see \cite{Su}, ch. 3, \S 2, (2.21).
\itemn{ii} If $p>0$, $2N-1=q$ is a power of $p$, $j=-1/4$ (including the case $N=3$, $p=5$, $j=-1/54$),
let $R\subset \bbF_{q^2}^\times$ be the subgroup of square roots of elements
of $\bbF_q^\times$, then
$${\rm Aut}_{k}(C)=\PGL_2(\bbF_q)\times_{\bbF_{\! q}^{^{\!\times}}} R$$
(the product is fibered with respect to the determinant
and the square map $R\to \bbF_{\! q}^{^{\!\times}}$).
\itemn{iii} In all other cases,
if $j\ne -1/4$ then ${\rm Aut}_{k}(C)=(\zmod{2})\times\bbD_N$,
and if $j=-1/4$ then ${\rm Aut}_{k}(C)=(\zmod{2})\times \bbD_{2N}$.
\end{trivlist}
\end{theo}

\begin{proo}
From now on we will denote by $O_\Gamma(x)$ the orbit of a point $x$ under
the action of a group $\Gamma$, and $\Gamma_x$ its stabilizer. We know that
$\tau$ has order 2 and is normal in $\Aut(C)$ (because of the
uniqueness of the hyperelliptic involution), hence it is central. Denote
$G=\Aut(C)/<\!\tau\!>$ so that there is a central extension
\begin{equation} \label{cent_ext}
1 \to <\tau> \to \Aut(C) \to G \to 1
\end{equation}
When $f\in \Aut(C)$, $f_0$ denotes its image in $G$. Denote by
$\Sigma=O_{\sigma_0}(\alpha)\cup O_{\sigma_0}(\beta)$ the set of the $2N$ branch
points of $\tau$, then by \ref{X/G_<=>_(Y,L,s)}, $G$ can be identified with the
subgroup of ${\rm Aut}(\bbP^1)$ of the homographies stabilizing $\Sigma$. What
we shall do is to determine $G$ thanks to Dickson's list, and then find the
class of the corresponding extension.

Notice that $(\zmod{2})\times\bbD_N\simeq <\!\tau,\sigma,\mu\!>\subset \Aut(C)$
so that $G$ contains a subgroup isomorphic to $\bbD_N$. Also $G$ acts
transitively on $\Sigma$, because $\bbD_N$ already does, and consequently
$$
\forall s \in \Sigma,\quad 2N=|O_G(s)|=|G|\,/\,|G_s|\;.
$$
Now, in four steps we read through the list in Dickson's theorem to find
all possible $G$'s. Before, we observe that we can allow conjugations of
$G$ in $\PGL_2(k)$ since it is just a change of variable on $x$ in equation
(\ref{EQ2}), so it does not change $C$ up to isomorphism. That is why from
step~2 on, we shall identify $G$ with the representatives given
in corollaries~\ref{cyclic_dihedral}, \ref{S4}, \ref{PGL2/PSL2}.

\titre{1st step: the groups that can not appear.}

	First of all, neither the cyclic groups nor $\frA_4$ possess dihedral
	subgroups $\bbD_N$, so they are ruled out.
	Now, let us see that the occurence of $G\simeq\frA_5$ is also impossible.
	The only dihedral subgroups in $G$ are $\bbD_3$ and $\bbD_5$, whence $N=3$ or 5.
	Assume $N=3$; then we must have $G_s\simeq \bbD_5$ for some $s\in\Sigma$, because
	it is a subgroup of $\frA_5$ with 10 elements. Denote $H:=<\!\sigma_0,\mu_0\!>
	\simeq \bbD_3$; then one can show that $G_s\cap H\simeq \zmod{2}$, and then
	there must be in $H$ a non trivial automorphism with a fixed point in $\Sigma$. This is
	a contradiction. When $N=5$, just interchange the roles of $\bbD_3$ and $\bbD_5$ (the
	former stands for $G_s$, the latter for $H$), and the argument carries on. Also, note
	that this is true for any $p$ (cases 1 and 3 of Dickson's theorem \ref{Dickson}).

	Finally, a group of type $Q\rtimes C$ (case 2) can not contain $\bbD_N$ with $(N,p)=1$.

\titre{2nd step: the case of $G=\PSL_2(\bbF_q)$ or $\PGL_2(\bbF_q)$, $q=p^s$.}

	We proceed in three substeps.
	\begin{trivlist}
	  \itemn{a} We show that $\Sigma$ coincides with the set of fixed points of
		order $p$ elements of $G$. Let $Q$ be a $p$-Sylow of the stabilizer
		$G_\alpha$, it is also a $p$-Sylow of $G$. If $g\in G$ has order $p$,
		then it belongs to some $p$-Sylow $Q'$. All Sylow subgroups are conjugated
		in $G$, so $Q'=tQt^{-1}$ and $g$ fixes $t(\alpha)\in\Sigma$. Conversely,
		if $s\in\Sigma$, then $G_s$ contains a $p$-Sylow
		of $G$, hence an element of order $p$. By cor.~\ref{fixed_points},
		$\Sigma=\bbP^1(\bbF_q)$; in particular, $2N=|\Sigma|=q+1$.
		In the sequel $\phi\in\bbF_{q^2}$ is a $2N$-th root of unity
		such that $\zeta=\phi^2$; observe that $\phi+\phi^{-1}$ and
		$\zeta+\zeta^{-1}$ both belong to $\bbF_q$.
	  \itemn{b} We work out the curve and the group structure of $\Aut(C)$.
		By (a) we get the birational equation $y^2=x^q-x$
		with genus $\frac{q-1}{2}=N-1$. Clearly $\PGL_2(\bbF_q)\subseteq G$ since
		it stabilizes $\bbP^1(\bbF_q)$. A priori $G$ could be bigger, however the
		only subgroups of $\PGL_2(k)$ containing $\PGL_2(\bbF_q)$ are
		isomorphic to $\PSL_2(\bbF_{q'})$ or $\PGL_2(\bbF_{q'})$, but then
		$q'=2N-1=q$ and hence $G=\PGL_2(\bbF_q)$. In particular the group
		$\PSL_2(\bbF_q)$ is ruled out. Moreover, for an element $\gamma=$
		\scriptsize
		$\left(
		\begin{array}{cc}
		a & b \\
		c & d
		\end{array}
		\right)
		$
		\normalsize
		and $\delta$ such that $\delta^2=\det(\gamma)$, there is an automorphism
		$$f_{\gamma,\delta}(x,y)= \left(\frac{ax+b}{cx+d}\:,
		\:\frac{\delta\:y}{(cx+d)^{\frac{q+1}{2}}}\right)$$
		This is the fiber product structure stated in the theorem. We can now
		check that this is indeed a Potts curve by giving an element of order $N$
		in $\PGL_2(\bbF_q)$. For this just consider the matrix $M_\zeta$
		of~\ref{remaMatrix}.
	  \itemn{c} At last we compute the invariant. In order to do this the quickest
		way is to notice that $\Sigma=\bbP^1(\bbF_q)$ is mapped to the $2N$-th roots
		of unity via the following transformation
		$$
		\gamma(x)=\frac{-x+1+\zeta^{-1}}{x-1-\zeta}
		$$
		Indeed, $u=\infty$ maps to -1 and for $u\in\bbF_q$ one has
		$$\gamma(u)^{2N}=\gamma(u)^{q+1}
		=\left(\frac{-u^q+1+\zeta^{-q}}{u^q-1-\zeta^q}\right)
		\left(\frac{-u+1+\zeta^{-1}}{u-1-\zeta}\right)
		=1$$ using that $u^q=u$ and $\zeta^q=\zeta^{-1}$.
		With the new variable $v=\gamma(x)$ we get an equation
		$$w^2=v^{q+1}-1$$
		and the invariant is $j=-1/4$.
	\end{trivlist}

\titre{3rd step: the case of $G={\sf S}_4=\:<ix,\frac{x+1}{x-1}>$ (for $N=3$,
see corollary~\ref{S4}).}

	This case can happen only if $p\ne 3$, cf theorem~\ref{Dickson}.
	Fix $\sqrt{i}$ a 8th root of unity and $i$ its square. With the notations of
	proposition~\ref{S4}, $x\mapsto ix$ is $\nu(1234)$ and $x\mapsto\frac{x+1}{x-1}$
	is $\nu(12)$. Applying a permutation on $\{1,2,3,4\}$ if necessary, we identify
	$\sigma_0$ and $\mu_0$ with $\nu(123)$ and $\nu(12)$, respectively. For $s\in\Sigma$,
	its stabilizer $G_s$ has cardinal 4, hence it is cyclic because two commuting
	involutions of $\bbP^1$ never have a common fixed point. We can choose $s$ so that
	$G_s$ is generated by $a_0(x)=ix$.

	Now $a_0$ can not have a single fixed point lying in $\Sigma$, because else
	it would act freely on $\Sigma-\{s\}$, and hence its order would divide $2N-1=5$.
	So the two fixed points $\{0,\infty\}$ are in $\Sigma$. The $G$-orbit of
	$\{0,\infty\}$ is $\{0,\infty,\pm 1,\pm i\}=\Sigma$. We can now give an equation
	$$y^2=x(x^4-1)$$
	with the automorphisms
	$$
	\begin{array}{cc}
	a(x,y)=(ix,\sqrt{i}\,y)
	& \sigma(x,y)=(-i\frac{x-1}{x+1},\frac{2\sqrt{2i}\,y}{(x+1)^3}) \\
	\mu(x,y)=(\frac{x+1}{x-1},\frac{2\sqrt{2}\,y}{(x-1)^3})
	& \lambda(x,y)=(-\frac{x-1}{x+1},\frac{2\sqrt{2}\,y}{(x+1)^3}) \\
	\end{array}
	$$
	After after a conjugation changing $\sigma_0$ into $x\mapsto jx$ ($j\in\mu_3^*$)
	we obtain the "Potts" equation $y^2=(x^3-(2+\sqrt{3})^3)(x^3+1)$. It allows to
	compute the invariant $j=-1/54$ but is less workable for the determination of
	the class of the extension~(\ref{cent_ext}). It can happen that
	$G\supsetneq {\sf S}_4$ only if $G$ is ${\rm PSL}_2(\bbF_q)$ or
	${\rm PGL}_2(\bbF_q)$, and we know that this implies $q=2N-1=5$. When $p=q=5$,
	we have $-1/54=-1/4$, and according to what was done before, $G=\PGL_2(\bbF_5)$.

	Now let us check, refering to the definition in \cite{Su}, chap.~2, \S~9, def.~9.10,
	that $\Aut(C)$ is a representation group of $\frS_4$. This is the
	last step to prove (i) in the theorem; as it is quite technical and not
	needed in the sequel we remain sketchy. We adopt the notations of~\cite{Su}.
	Denote $H=\Aut(C)$ and $Z=<\tau>$. Assume $L$ is a proper subgroup of $H$
	such that $H=ZL$, then $[H:L]=2$, $L\cap Z=1$, hence $H=Z\times L$ does not
	contain any element of order 8, contradicting the existence of $a$. Furthermore,
	denote by $M(G)$ the Schur multiplier, then $|Z|=|H'\cap Z|=2=|M(G)|$ by
	checking that $\tau=[\mu,\lambda]$ is a commutator. Third, obviously,
	$|H|=|G|.|M(G)|$. Then the result follows by \cite{Su}, chap. 3, \S 2, (2.21).

\titre{4th step: the dihedral case $G={\sf D}_M=<\varepsilon x,\frac{1}{x}>$,
$\varepsilon\in\mu_M^*$, for $N|M$, see corollary~\ref{cyclic_dihedral}.}

	It is the last possibility. For $w\in\bbP^1$ its ${\sf D}_M$-orbit is
	$\left\{w,\varepsilon w,\dots,\varepsilon^{M-1}w,\,\frac{1}{w},
	\dots,\frac{\varepsilon^{M-1}}{w}\right\}$. This
	has cardinal 2 if $w=0$ or $\infty$, $2M$ if $w$ is in general position,
	and $M$ if $w$ is one of $\pm\sqrt{\varepsilon^j}$. We conclude that, in
	general, $M=N$; $M=2N$ occurs when the points of the orbit are vertices
	of a regular $2N$-gon, yielding the equation $y^2=x^{2N}-1$. The invariant
	is then $j=-1/4$, and $\sigma_0$ has a "root" $\sqrt{\sigma_0}(x)=\zeta_{_{2N}} x$.
\end{proo}

\begin{coro} \label{coroAutMod}
For all $N,p$, the group of automorphisms of $(C,\sigma)$ is:
\begin{center}\begin{tabular}{rl}
$j\ne -1/4\,:$ & ${\rm Aut}_{k}(C,\sigma)=(\zmod{2})\times (\zmod{N})$ \\
$j=-1/4\,:$ & ${\rm Aut}_{k}(C,\sigma)=(\zmod{2})\times (\zmod{2N})$.
\end{tabular}\end{center}
\end{coro}

\begin{proo}
The only non-obvious case is (ii) of the theorem. We keep the notations
of step 2 in the proof of the theorem. It suffices to check
that the centralizer $Z$ of $\sigma_0$ in $\PGL_2(\bbF_q)$ is cyclic of order
$2N$. It is cyclic because, after a conjugation changing $\sigma_0$
into $x\mapsto \zeta x$, $Z$ maps to a finite subgroup of $\bbG_m$.
It has order $\le 2N=q+1$ by corollary~\ref{orders_in_PGL2(Fq)},
observing that $\sigma_0$ is a power of a generator of $Z$. Finally it has
order exactly $2N$ because there is in $\PGL_2(\bbF_q)$ an explicit square
root for $\sigma_0$. Indeed, if $\sigma_0$ is given as above by the matrix
$M_\zeta$ of~\ref{remaMatrix}, whose Cayley-Hamilton polynomial is
$X^2-(\zeta+\zeta^{-1})X+(\zeta+\zeta^{-1})$, one finds that
$M_\zeta+(\phi+\phi^{-1})\Id$ has a square equal to $M_\zeta$ in $\PGL_2(\bbF_q)$.
\end{proo}


\subsection{Wild case} \label{CS}

Now we study $N$-Potts curves when $p|N$. The ground field $k$ is
still assumed to be algebraically closed, of characteristic $p>0$.
Here it is definition~\ref{def0} that applies: the ramification data of
definition~\ref{def_Potts_curve} do not make sense any more. As
usual $\sigma$ stands for the given automorphism of order $N$ and
$\tau$ is the hyperelliptic involution. Observe
that if we have an isomorphism $\varphi:C\to C'$ between two Potts
curves with $\varphi\sigma=\sigma'\varphi$, then it follows from unicity
of the hyperelliptic involution that we also have $\varphi\tau=\tau'\varphi$.

As a matter of fact it is not so clear that such curves exist. Let
$(C,\sigma,\tau)$ be an $N$-Potts curve with arbitrary $N$ multiple of $p$.
Then $\sigma$ induces an automorphism of order $N$ on the quotient
$C/\tau\simeq\bbP^1$. By proposition~\ref{automorphismes_PGL2},
the only possible case is $$N=p\;,$$ and the induced automorphism
$\sigma_0$ is $x\mapsto x+1$ up to conjugation.
As in proposition~\ref{Potts=>hyperelliptique}, it is easy to count $2p$ fixed
points for the hyperelliptic involution $\tau$, they form two orbits of $\sigma$
(just lift the 2 fixed points of the involution induced on $C/\sigma$).
The corresponding affine model is

\begin{equation} \label{EQ4}
y^2=(x^p-x)^2+A(x^p-x)+B \qquad \qquad (A,B \in k).
\end{equation}

\noindent Hence, there exist $N$-Potts curves with $p|N$ exactly when $N=p$.
In that case we define a modular invariant by
$$
j=\frac{1}{A^2-4B}
$$

\begin{prop} \label{`couple'_de_Potts}
Two $p$-Potts curves $(C,\sigma,\tau)$ and $(C',\sigma',\tau')$ of invariants
$j$ and $j'$ are isomorphic if and only if $j=j'$.
\end{prop}

\begin{proo}
Let $\varphi:C\to C'$ be an isomorphism with $\sigma'\varphi=\varphi\sigma$.
It induces a map $\tilde{\varphi}:C/\tau\to C'/\tau'$ on the quotients.
Moreover, again by~\ref{automorphismes_PGL2} we can assume that both $\sigma$
and $\sigma'$ are $x\mapsto x+1$. So if we consider equations of type~(\ref{EQ4})
for $C$ and $C'$, then $G$-equivariance reads $\tilde{\varphi}(x+1)=\tilde{\varphi}(x)+1$.
From this we deduce that $A'=A+2(t^p-t)$ and $B'=B+(t^p-t)A+(t^p-t)^2$,
and then $j'=j$.

Conversely if $j'=j$ then the choice of a root $t\in k$ of $t^p-t=\frac{1}{2}(A'-A)$
satisfies also $B'=B+(t^p-t)A+(t^p-t)^2$. This ensures that $\varphi(x,y)=(x+t,y)$
defines an equivariant isomorphism between $C$ and $C'$ with equations~(\ref{EQ4}).
\end{proo}

\begin{rema}
Here, contrary to the tame case (compare with~\ref{j-isom}) there is no
numerical invariant such as $[\chi]$ but only a continuous one. The computation
of the moduli space in both cases will enlighten this in the next sections.
\end{rema}

At last we compute the automorphism group of $p$-Potts curves.
Let $(C,\sigma,\tau)$ be given by equation~(\ref{EQ4}). Let $r,s$ be
the roots of $T^2+AT+B$, and $\alpha$ (resp. $\beta$) a root of $T^p-T-r$
(resp. $T^p-T-s$), so that~(\ref{EQ4}) reads
$$
y^2=(x^p-x-r)(x^p-x-s)=\prod_{i=0}^{p-1}\,(x-\alpha+i)\:\prod_{i=0}^{p-1}\,(x-\beta+i)\;.
$$
We have the following automorphisms:
$$
\begin{array}{ccc}
\sigma(x,y)=(x+1,y) \quad; &
\tau(x,y)=(x,-y) \quad; &
\mu(x,y)=(\alpha+\beta-x,y) \;.
\end{array}
$$

\begin{prop} \label{automorphismes_car_p}
Let $(C,\sigma,\tau)$ be a $p$-Potts curve, then
$\Aut_k(C)\simeq (\zmod{2})\times\bbD_p\;$ and
$\Aut_k(C,\sigma,\tau)\simeq (\zmod{2})\times (\zmod{p})\;$.
\end{prop}

\begin{proo}
It is still true (cf theorem~\ref{theo_autom}) that
\begin{trivlist}
\itemm{i} $\tau$ is of order 2, normal and central,
\itemm{ii} $G=\Aut(C)/<\!\tau\!>$ is the subgroup of $\Aut(\bbP^1)$
of homographies stabilizing $\Sigma=O_{\sigma_0}(\alpha)\cup O_{\sigma_0}(\beta)$,
\itemm{iii} $\bbD_p\subset G$ and $2p=[G:G_x],\quad\forall x\in\Sigma$.
\end{trivlist}
Let $Q$ be a $p$-Sylow of $G$ containing $\sigma_0$; by Dickson's theorem $Q$
is elementary abelian. Assume that $Q$ has order more than $p$, then there exists $\theta\in Q$
commuting with $\sigma_0$, hence $\theta(x)=x+u$, with $u\not\in \bbF_p$. Moreover
$\theta$ stabilizes $\Sigma$, and exchanges the orbits $O_{\sigma_0}(\alpha)$ and
$O_{\sigma_0}(\beta)$ since $u\not\in\bbF_p$. Therefore
$$
(\exists i,j\in\bbF_p)\qquad
\left\{
\begin{array}{l}
\alpha+u=\beta+i \\
\beta+u=\alpha+j 
\end{array}
\right.
$$
which implies $u=\frac{i+j}{2}\,\in\bbF_p$, a contradiction. Consequently
$Q=<\sigma_0>$; we can now read through the list in Dickson's theorem.

If $G=\PSL_2(\bbF_p)$ or $\PGL_2(\bbF_p)$, then $G_\alpha$ has order $(p^2-1)/d$
with $d=2$ or 4. This is prime to $p$, hence $G_{\alpha}$ est cyclic,
but this contradicts corollary~\ref{orders_in_PGL2(Fq)}.

The only remaining possibility is $G=Q\rtimes C=<\sigma_0>\rtimes <\varphi>$,
because $\frA_5$ is ruled out by the same arguments as in the case
$(N,p)=1$. The order of $\varphi$, denoted $2m$, is prime to $p$; changing the point
$\alpha$ in the orbit if necessary, we can assume that its stabilizer is
$G_{\alpha}=<\varphi^2>$. As $Q$ is normal, $\varphi^{-1}\sigma_0\varphi=\sigma_0^\ell$
for some $\ell\in\bbF_p^*$, from which we deduce that $\varphi(x)=\ell^{-1}x+b$
for some $b\in k$. As $\varphi^2$ fixes
$\alpha$, we derive $(\ell^2-1)\alpha=\ell(\ell+1)b$. If $\ell+1\ne 0$ it implies
$\varphi(\alpha)=\ell^{-1}\alpha+b=\alpha$; so $\ell=-1$, $\varphi^2=1$ and $m=1$.
\end{proo}

We see that the remarkable symmetry previously obtained when $j=-1/4$ does not
occur here in characteristic~$p$. In particular, this implies that the $p$-Potts
curve in characteristic~$0$ with invariant $j=-1/4$ can not have good reduction
in characteristic~$p$, i.e. that any model of $C$ over a discrete valuation
ring of residue characteristic $p$ will have a singular special fibre.


\section{The stack $\clP_N$ when $N$ is composite}
\label{construction_EM}

Let $N\ge 3$ be a non-prime integer. Let $k$ be an algebraically closed field of
characteristic $p\ne 2$. We established in the previous section that whenever $p$ is
prime to $N$, there is a bijection between isomorphism classes of $N$-Potts
curves and a sum of $\varphi(N)/2$ copies of the affine punctured line
${\bbA}^1_*:={\bbA}^1-\{0\}$ over $k$ (see~\ref{remaA1*}). Furthermore when
$p$ divides $N$ we saw in~\ref{CS} that there are no Potts curves at all. We now
show (th.~\ref{espace_modulaire}) that the coarse moduli space of the stack $\clP_N$ is
${\bbA}^1_*\otimes\bbZ[\frac{\zeta+\zeta^{-1}}{2},\frac{1}{2N}]$ (which indeed
splits as a disjoint sum over any field containing the $N$-th roots of unity).
Here a couple of remarks are in order:
\begin{trivlist}
\itemn{i} To be more precise, in all what follows, whenever we write $\bbZ[1/2,\zeta]$
we mean the ring of cyclotomic integers, and whenever we write
$\bbZ[1/2,\frac{\zeta+\zeta^{-1}}{2}]$, we mean its ring of invariants under
$\zeta\mapsto \zeta^{-1}$. In other words, the former ring is $\bbZ[1/2][X]/\Phi_N$
and the latter is $\bbZ[1/2][X]/\psi_N$ where $\Phi_N$ and $\psi_N$ are the
cyclotomic polynomials of~\ref{remaMatrix}.
\itemn{ii} It will become clear while reading that everything in this section
applies equally well to the tame stack $\clP_p\otimes\bbZ[1/2p]$ when $N=p$ is prime.
\end{trivlist}
As an immediate consequence we have a result of good reduction (\ref{reduction}).
In~\ref{GPP}, we compute the modular Picard group of the fibres of the stack $\clP_N$.
At last we determine topologically the stable curves that are involved as stable
limits in the process of compactification of the moduli space of $\clP_N$, in~\ref{CCS}.

\subsection{Preliminaries}

\begin{defi} \label{def_famille}
Let $S$ be a scheme over $\bbZ[1/2]$. An $N$-Potts curve over $S$ is a
triple $(C,\sigma,\tau)$ composed of a smooth projective $S$-curve, and
two automorphisms, $\sigma:C\to C$ of order $p$ and $\tau:C\to C$ of order 2,
such that the geometric fibers are Potts curves in the sense of definition~\ref{def0}.
\end{defi}

It follows from standard arguments that the stack $\clP_N$ of $N$-Potts
curves is a separated Deligne-Mumford stack over $\bbZ[1/2]$ (see~\cite{We}
and \cite{BR}). Also, as we saw in section~\ref{section_CP_autom},
when $N$ is invertible in the structure sheaves of the base schemes,
we might as well use :

\begin{defi}
Let $S$ be a scheme over $\bbZ[1/2N]$. An $N$-Potts curve
is a proper smooth morphism of schemes $f:C\to S$, together with an
$S$-automorphism $\sigma:C\to C$ of order $N$, such that the geometric
fibers $(C_s,\sigma_s)$ are Potts curves in the sense of
definition~\ref{def_Potts_curve}.
\end{defi}

This gives a much better understanding of the stack, so we will work with
this definition. Then let $f:C\to S$ be a Potts curve over $S$ (remark:
the type $[\chi]$ is locally constant over $S$). We will need to recover
the following important fact that $f:C\to S$ is a hyperelliptic curve :

\begin{lemm} \label{tau}
There exists an involution $\tau:C\to C$ such that $f:C\to S$
becomes a family of hyperelliptic curves (in the sense of {\rm \cite{KL}}).
\end{lemm}

\begin{proo}
Let $G=<\!\sigma\!>$. As $C$ is projective over $S$, the quotient $D=C/G$
exists; by smoothness the quotient map $\pi:C\to D$ is finite flat of degree $N$.
As $N\in \clO_S^\times$ its formation commutes with base change
(see~\cite{KaMa}, A7.1.3.4), as is also the case for the branch locus $B=\pi_*(C^G)$.
In particular $D$ is a $\bbP^1$-bundle over $S$ and $B$ is \'etale and finite
of degree 4 over~$S$. We may localize (in the \'etale topology) as much as
desired, because by unicity the hyperelliptic involution will exist globally.

Hence we may assume that $S=\Spec(R)$ is affine, that $D=\bbP^1_R$,
and $B$ is a sum of four disjoint sections $\alpha,\beta,\gamma,\delta$.
Write the sections as $\alpha=(a_{_1}:a_{_2}),\dots,\delta=(d_{_1}:d_{_2})$.
We define a linear transformation of $D$ by the matrix
$$
\tau_0=\left[
\begin{array}{cc}
a_{_1}b_{_1}c_{_2}d_{_2}-c_{_1}d_{_1}a_{_2}b_{_2} &
a_{_1}c_{_1}d_{_1}b_{_2}+b_{_1}c_{_1}d_{_1}a_{_2}
-a_{_1}b_{_1}c_{_1}d_{_2}-a_{_1}b_{_1}d_{_1}c_{_2} \\
(a_{_1}b_{_2}+a_{_2}b_{_1})c_{_2}d_{_2}-(c_{_1}d_{_2}+c_{_2}d_{_1})a_{_2}b_{_2}
& -(a_{_1}b_{_1}c_{_2}d_{_2}-c_{_1}d_{_1}a_{_2}b_{_2})
\end{array}
\right]
$$
(we mimic the expression in the case where the base is a field).
Its determinant
$\det(\tau_0)=-(a_{_1}c_{_2}-a_{_2}c_{_1})(a_{_1}d_{_2}-a_{_2}d_{_1})
(b_{_1}c_{_2}-b_{_2}c_{_1})(b_{_1}d_{_2}-b_{_2}d_{_1})$
is indeed invertible, as is clear fibrewise, and $\tau_0$ is involutive
because of the vanishing of the trace (cf~\ref{examsOrders}).
We must now lift it to~$C$. By~\ref{X/G_<=>_(Y,L,s)} the cyclic
covering $\pi:C\to D=\bbP^1_R$ is described by $\clL\simeq {\cal O}(-2)$
and a global section $s$ of $\clL^{ N}$, with
$\clL^{ N}\simeq {\cal O}(-{\cal D})$ where
${\cal D}=\alpha+\beta+(N-1)\gamma+(N-1)\delta$. By construction
$\tau_0$ respects the data $(\clL,s)$, hence it lifts
to an automorphism $\tau$ of $C$ with $\tau\sigma=\sigma\tau$. As in the case
of a single Potts curve, we can assume $\tau^2=\Id$. This completes the proof.
\end{proo}

\subsection{The moduli space} \label{dem_theoreme}

We now compute the moduli space of $\clP_N$. The proof below will in fact give
a concrete expression of $\clP_N$ as a quotient stack of an open subset in the
affine 3-space of homogeneous polynomials of the form $H(X,Z)=UX^{2N}+AX^NZ^N+BZ^{2N}$
by the action of the group $(\bbG_m)^2$. The first $\bbG_m$ factor acts simply
(and freely) by multiplication, and the second factor acts by
$\lambda.H(X,Z):=H(\lambda X,Z)$ (this is just the isomorphism relation between Potts curves,
see proposition~\ref{j-isom}). This description is at least totally correct over
an algebraically closed. Being rather interested in the arithmetic and reduction
of the moduli space, we will not insist on this aspect. Hence we will show~:

\begin{theo} \label{espace_modulaire}
The moduli space of $\clP_N$ is
${\bbA}^1_*\otimes\bbZ[\frac{\zeta+\zeta^{-1}}{2},\frac{1}{2N}]$.
\end{theo}

\noindent To avoid heavy notations we will write $P$ for
${\bbA}^1_*\otimes\bbZ[\frac{\zeta+\zeta^{-1}}{2},\frac{1}{2N}]$.
We split the proof into two steps~:

\titre{1st step: we build the morphism to the moduli space.}

Here again we will work \'etale locally on $S$, and take care that the
construction of the morphism is canonical enough so that it descends.
Notations are as above. By~\ref{tau} there is an
involution $\tau\in{\rm Aut}_S(C)$. Then $E=C/\tau$ is a
$\bbP^1$-bundle over $S$, we denote by $r:C\to E$ the natural
projection and by $q:E\to S$ the structure morphism. As $\sigma$
commutes with $\tau$ it induces an automorphism of order $N$ on $E$.
The divisor of fixed points $T=E^\sigma$ is an \'etale cover of degree~2
of $S$, \'etale locally it is a sum of two disjoint sections
$\Delta+\Delta'$. Choosing one of the two (say $\Delta$)
defines an invertible sheaf $\clO(\Delta)$ of degree 1. We may call
it $\clO(1)$ and then $E\simeq \bbP(V)$ with $V=q_*\clO(1)$. By
disjointness, if $L$ and $M$ are the restrictions of
$\clO(1)$ to $\Delta$ and $\Delta'$ respectively (viewed as
sheaves on $S$), then $V$ splits as $L\oplus M$.
By construction the action of $\sigma$ is now diagonal, given by
multiplication by two invertible global sections $s,t\in\Gamma(S,\clO_S)^\times$.
We can normalize by changing $V$ into $V\otimes L^{-1}$, so that $L=\clO_S$
and $s=1$. Then as $\sigma$ has order $N$, $t$ is a primitive $N$-th root of unity.

Now let us consider the double cover $r:C\to E$. It is described by the
decomposition $r_*{\cal O}_C={\cal O}_E\oplus \clL$ and by a "Weierstrass"
section $\theta\in\Gamma(E,\clL^{ -2})$, well determined up to an element of
$\Gamma(E,\clO_E^\times)$ (see~\ref{X/G_<=>_(Y,L,s)}). Since $\clL^{ -1}$ has
degree~$N$ on the fibres, we have $\clL^{ -1}\simeq \clO(N)\otimes q^*K$ for
some $K\in \Pic(S)$. Also as $\theta_{|\Delta}$ is everywhere nonzero (the
fixed loci for the actions of $G$ and $\tau$ are disjoint fibrewise), by restriction
to $\Delta$ we get $K^{ 2}\simeq \clO_S$, hence $\clL^{ -2}\simeq \clO_E(2N)$.
Consequently we can identify $\theta$ with a $\sigma$-invariant section of
$$
 \Gamma(E,{\cal O}_E(2N))
=\Gamma(S,{\rm Sym}^{2N}(V))
=\bigoplus_{j=0}^{2N}\,\Gamma(S,M^{ j})\;.
$$
The most convenient writing is to use global coordinates
$X,Z$ on $\bbP(V)$, with say $X$ corresponding to $M$. Then
$\theta\in\Gamma(S,{\cal O}_S\oplus M^{ N}\oplus M^{ 2N})$, so
$\theta\propto H=UX^{2N}+AX^NZ^N+BZ^{2N}$ for some sections $U,A,B$ of $\clO_S$,
$M^N$ and $M^{2N}$ respectively (recall that $\propto$ means equality up to an
invertible element). Looking locally on the fibers, one sees that
$U$ is invertible. Finally note that the smoothness of the fibres of $C/S$
requires that the section $\theta$ has no multiple zero (on all fibres).
This means that neither $B$ nor $A^2-4UB$ vanish, hence there is a
well-defined section $j=\frac{UB}{A^2-4UB}\in\Gamma(S,{\cal O}_S)^\times$.
The assignments $F(\zeta)=t$ and $F(X)=j$ give a map
$F:\bbZ[\zeta,\frac{1}{2N}][X,X^{-1}]\to \Gamma(S,{\cal O}_S)^\times$.

Now we proceed to check the independance of this map with respect
to the choices made. In fact the only place where there is a different
possibility is the choice of $\Delta$ rather than $\Delta'$.
Choosing $\Delta'$ is equivalent to exchanging the
coordinates $X,Z$ on $\bbP(V)$, so clearly $j$ is unchanged. However
$t$ is changed into $t^{-1}$; so in any case if we set
$F(\zeta+\zeta^{-1})=t+t^{-1}$ then the map
$$
\bbZ[\frac{\zeta+\zeta^{-1}}{2},\frac{1}{2N}]
\left[X,X^{-1}\right]\to \Gamma(S,{\cal O}_S)^\times
$$
is independant of the choice. Eventually we have a map
$S\to {\bbA}^1_*\otimes\bbZ[\frac{\zeta+\zeta^{-1}}{2},\frac{1}{2N}]=P$.
It is clear that the construction is functorial, providing a morphism
$\Phi:\clP_N\to P$.

\titre{2nd step: we check the properties of the moduli space.}

Recall that we must check two things: first, that for any geometric point
$\Spec(k)\to \Spec(\bbZ[1/2])$, the morphism $\Phi$ induces a bijection
between isomorphism classes in $\clP_N(k)$ and $P(k)$. If the characteristic
of~$k$ divides~$N$ then it is clear because both are empty, and else this is
exactly~\ref{remaA1*} (observe that $P\otimes k$ splits as $\varphi(N)/2$
copies of $\bbA^1_*(k)$). The second thing is to see that every map from
$\clP_N$ to an algebraic space factors through $P$. This can be done as
in~\cite{MS}, using a family whose classifying morphism $S\to P$ is finite
surjective (such a family is sometimes called \emph{tautological}). For
this we consider the one-parameter family $C_0$ with equation
$y^2=x^{2N}+\lambda x^N+1$, over the base
$S_0=\Spec(\bbZ[\zeta,\frac{1}{2N}][\lambda,\frac{1}{\lambda^2-4}])$.
(This is of course only an affine smooth curve; to make the definition
rigorous we actually glue $C_0$ with another copy of itself, compatibly with
the maps to $S_0$, along the open set $x\ne 0$, via the isomorphism $x'=1/x$,
$y'=y/x^N$). The data of $\zeta$ and of the invariant $j_0=\frac{1}{\lambda^2-4}$
determine a morphism $S_0\to P$ which we denote by $\Phi_0$. Now let
$\Psi:\clP_N\to Q$ be a morphism of stacks to an algebraic space, and let
$\Psi_0:S_0\to Q$ be the morphism corresponding to $\Psi(C_0)$; let
$\Gamma\subset P\times Q$ be the scheme-theoretic image of $h=(\Phi_0,\Psi_0)$.

\begin{diagram}[height=1.6em,width=2em,PostScript=dvips]
& & S_0 & & \\
& \ldTo(2,4)^{\Phi_0} & \dTo^{_h} & \rdTo(2,4)^{\Psi_0} & \\
& & \Gamma & & \\
& \ldTo_{p_1} & & \rdTo_{p_2} & \\
P & & \rDashto & & Q
\end{diagram}

\noindent We observe that, $\Phi_0$ being finite and $p_1$ separated, $h$ is finite.
In particular, $h$ is closed, hence $\Gamma=h(S_0)$ as sets. Second notice that
$\Gamma$ is integral because $S_0$ is.
Third $p_1$ is closed and bijective: it is closed and surjective because
$\Phi_0$ is, and injective because for $s,s'\in S_0$,
$\Phi_0(s)=\Phi_0(s') \Rightarrow C_{0,s} \simeq C_{0,s'}
\Rightarrow \Psi_0(s)=\Psi_0(s')$ ($Q$ being a space).
Thus $p_1$ is dominant, bijective and separable
(since $\Phi_0$ is), hence it is a birational map. At last, as $P$ is normal,
Zariski's Main Theorem states that $p_1$ is an isomorphism. Then the composition
$p_2\circ p_1^{-1}:P\to Q$ gives a morphism
which factors~$\Psi$. \qed

Of course it must be said that $P$ is not a fine moduli space. This is due to
the presence of automorphisms; actually, in view of~\ref{coroAutMod} we could
get rid of the group $(\zmod{2})\times (\zmod{N})$ by a process to be explained
in~\ref{Twoquotient}, but the extra automorphism when $j=-1/4$ still causes
ramification of $j$ above $-1/4$:
$$
j+\frac{1}{4}=\frac{A^2}{4(A^2-4B)}
$$

A consequence of the explicit construction of the classifying morphism
is that obviously reduction modulo $\frp$ can be done at any prime $\frp>2$.
The result is straightforward:

\begin{theo} \label{reduction}
Assume that $N\ge 3$ is a composite integer and that $\frp>2$ is a prime.
Then the moduli space of $\clP_N$ has good reduction at $\frp$, i.e. the
moduli space of $\clP_N\otimes \bbF_\frp$ is the (possibly empty) fibre
at $\bbF_\frp$ of the moduli space of $\clP_N$. \qed
\end{theo}

\subsection{The Picard group} \label{GPP}

Let $k$ be an algebraically closed field of characteristic $p$ prime to $2N$.
We are now going to compute the Picard group of the geometric fibre 
$\clP_N\otimes k$. This will, in some sense, reveal that the \emph{geometry}
of the stack carries the dependance on $N$ (whereas the subring of cyclotomic
integers involved in the moduli space is of an \emph{arithmetic} nature).
Actually, as $\clP_N\otimes k$ splits as a sum of isomorphic stacks
$(\clP_N\otimes k)_{[\chi]}=\clP_{N,[\chi]}$ according to the classes
of characters $[\chi]$, we will compute the Picard group of one of them.

We first briefly recall the definitions. The \emph{\'etale site} of $\clP_N$
has as open sets the \'etale morphisms $u:U\to \clP_N$ from an algebraic space,
denoted $(U,u)$ or $U$. A map between two open sets $U,V$ is a couple $(f,\alpha)$
with a morphism of algebraic spaces $f:U\to V$ and a 2-isomorphism
$\alpha:v\circ f\isomto u$. Equivalently, it is a 2-commutative triangle
with vertices $U,V,\clP_N$. Briefly said, the coverings of $(U,u)$ are
the \'etale, surjective families $\coprod U_i\to U$, and the topology generated
by all these is the \emph{\'etale site} $(\clP_N)_\et$. Finally an
\emph{invertible sheaf} $L$ on $(\clP_N)_\et$ is given by a collection
of invertible sheaves $L|_U$ on $U$ for every open set $U\to \clP_N$, and
isomorphisms $\ell_{f,\alpha}:f^*L|_V\isomto L|_U$ for all maps as above
between open sets, such that any composition
$$
U\stackrel{(f,\alpha)}{\longrightarrow}V\stackrel{(g,\beta)}{\longrightarrow}W
$$
gives rise to an equality
$\ell_{g\circ f,\alpha\circ f^*\beta}=\ell_{f,\alpha}\circ f^*\ell_{g,\beta}$.
The Picard group is the set of isomorphism classes of invertible
sheaves, endowed with the obvious tensor product.

\begin{lemm} \label{LemmBeta}
There is a morphism of groups
$\beta:{\rm Pic}(\clP_{N,[\chi]})\to \zmod{2}\times \zmod{2N}$.
\end{lemm}

\begin{proo}
We define it as Mumford does in~\cite{Mu}. Let $L$ be an invertible sheaf on
$\clP_{N,[\chi]}$ and $U\to \clP_{N,[\chi]}$
an open set, corresponding to an $N$-Potts curve $(C,\sigma,\tau)$ over $U$
(we use definition~\ref{def0}). In terms of the topology on $\clP_{N,[\chi]}$,
$\sigma$ gives an automorphism of the open set $U$,
so that there is an isomorphism $\ell_{\Id,\sigma}:L|_U\isomto L|_U$. It is
given by an invertible global section of $\clO_U$, and actually by the
compatibility of $L$ w.r.t composition this section is an $N$-th root of
unity. Finally we get a morphism $U\to\mu_N$ to the scheme $\mu_N$ of
$N$-th roots of unity.

It is clear that for a connected $U$ the image in $\mu_N$ is constant;
it is even independant of the chosen open set $U$, because the stack
$\clP_{N,[\chi]}$ is irreducible, so two nonempty open sets $(U,u)$
and $(V,v)$ have images that intersect in $\clP_{N,[\chi]}$.
In particular, if we choose $U$ to be an atlas, one $k$-point $x\in U$
will give the curve $C_x$ with the extra automorphism $\sigma_0$ whose
square is $\sigma$ (for the value $-1/4$ of the invariant). As above,
we then get a point in $\mu_{2N}$, and obviously its square is the
point in $\mu_N$ computed before.

Now recall that a primitive $N$-th (resp. $2N$-th) root of unity $\zeta$
(resp. $\phi=-\zeta$) is determined up to inversion by $[\chi]$. This
yields an isomorphism $\mu_{2N}\simeq \zmod{2N}$ mapping $\phi$ to 1,
hence for given $L$ there is a well-defined element $\beta_2(L)\in\zmod{2N}$
computed with any nonempty open set~$U$. The same works with the
hyperelliptic involution $\tau$; in order to keep additive notation
we define $\epsilon=\beta_1(L)\in\zmod{2}$ to be such that
$\ell_{\Id,\tau}$ is the multiplication by $(-1)^\epsilon$. This
completes the definition of $\beta=(\beta_1,\beta_2)$.
\end{proo}

As noticed in the proof, $\beta$ is determined up to inversion of the
generator in $\zmod{2N}$. The following should now be quite close to intuition:

\begin{lemm} \label{beta_surjectif}
$\beta$ is surjective.
\end{lemm}

\begin{proo}
Let $U\to \clP_{N,[\chi]}$ be an atlas, and $(f:C\to U,\tau,\sigma)$
be the corresponding curve. It is tempting, as in \cite{Mu}, to evaluate
$\beta$ on the Hodge bundle $L=\bigwedge^{N-1}f_*\Omega_C$ where
$\Omega_C$ is the sheaf of differential 1-forms. Clearly we can compute
$\beta$ on the fibre of $L$ over a single point of $U$, and of course we
choose a point $x\in U$ whose fibre is the curve $C_x$ with an extra
automorphism $\sigma_0$ (for $j=-1/4$). Up to isomorphism $C_x$ has equation
$y^2=x^{2N}-1$. The basis of $\Gamma(C,\Omega_C)$ given by the forms
$\omega_i=x^{i-1}\frac{dx}{y}$ for $1\le i\le N-1$, has the virtue of
diagonalizing the action of the group $G\simeq\zmod{2}\times \zmod{2N}$
generated by $\tau$ and $\sigma_0$. Indeed
$$
\begin{array}{rcr}
\tau(\omega_i) & = & -\omega_i \\
\sigma_0(\omega_i) & = & \phi^i\omega_i
\end{array}
$$
Unfortunately we see that $\sigma_0$ maps $\omega_1\wedge\dots\wedge\omega_{N-1}$
to $(-1)^{(N-1)/2}$ times itself (as does $\tau$) so we can not conclude.
Thus taking maximal exterior power was too crude, but instead we can use
the fact that $f_*\Omega_C$ is a $G$-sheaf, so
$$
f_*\Omega_C=\oplus_{i=1}^{N-1}\,L_i
$$
where $L_i={\rm ker}(\sigma^*-\phi^i\,{\rm id})$ is an invertible sheaf.
We obtain $\beta(L_1)=(1,1)$ and $\beta(L_2)=(1,2)$, that generate the image.

\end{proo}

The end result is very similar to the one in the elliptic case (see~\cite{Mu}):

\begin{theo} \label{PicardTame}
$\beta$ is injective, i.e.
$\Pic((\clP_N\otimes k)_{[\chi]})\simeq (\zmod{2})\times (\zmod{2N})$.
\end{theo}

\begin{proo}
Given an invertible sheaf $L$ such that $\beta(L)=0$, we show that $L|_U$
is trivial for any $u:U\to \clP_N$. Let $f:C\to U$ be the corresponding
curve, and let us simplify the notations to $\clP:=\clP_{N,[\chi]}$ and $L:=L|_U$.
As in~\cite{Mu} we will show that $L$ "descends" to the moduli space $P$, but
we will write down carefully the argument (only allusive in~\cite{Mu}) since
it involves nonflat descent. Consider the diagram
\begin{diagram}[height=1.6em,width=2em,PostScript=dvips]
U\times_{\clP} U & \rTo^a &
U\times_{P} U & \rInto^b &
U\times_S U & \pile{\rTo^{p_1}\\ \rTo_{p_2}} & U
\end{diagram}
Let $q_i:=p_i\circ b\circ a$, then by definition of an invertible sheaf
we have an isomorphism $\ell:q_1^*L\isomto q_2^*L$ between sheaves
on $U\times_{\clP} U$. But the latter space is just
$I:=\Isom(p_1^*C,p_2^*C)$, finite and unramified over $U\times_S U$.
By definition of the coarse moduli space, $U\times_{P} U$
is the image of $I$ in $U\times_S U$. So actually $I$ is a "torsor"
over $U\times_{P} U$, with structure group $G=\Aut((p_1\circ b)^*C)$.
We must be careful that $G$ is not flat, so the meaning of a torsor
here is just that $I\times I\simeq G\times I$. We have as usual an invariant
pushforward $a_*^G$ (pushing forward and then taking invariant sections)
and it satisfies $a_*^Ga^*F\simeq F$ for any locally free sheaf $F$ of finite
rank on $U\times_{P} U$. Here flatness of $G$ is indeed unnecessary,
if $F$ is locally free: using that both $I$ and $G$ are affine over
$U\times_{P} U$, we may locally reduce to the following situation of
commutative algebra. We have $U\times_{P} U=\Spec(A)$, $I=\Spec(B)$,
$G=\Spec(A[G])$ with a coaction $B\to A[G]\otimes B$ such that $A=B^G$;
finally $F$ is given by a free module $M$. It remains to check that
$(M\otimes_A B)^G=M$ which is clear since $M$ is free.

The initial assumption that $\beta(L)=0$ says exactly that
$\ell:q_1^*L\isomto q_2^*L$ is a $G$-equivariant isomorphism,
so applying $a_*^G$ we obtain an isomorphism
$\psi:(p_1\circ b)^*L\isomto (p_2\circ b)^*L$. Fortunately, the stack
$\clP$ as well as its moduli space $P$ are smooth,
therefore the map $\clP\to P$ is flat. Thus the
map $U\to P$ is flat, and $\psi$ is a descent datum for $L$.
Finally $L$ descends to $P$, so it is trivial since
$\Pic(P)=\Pic(\bbA^1_*\otimes k)=0$.
\end{proo}

\subsection{Compactification by stable curves} \label{CCS}

It is known by the general theory of tame Hurwitz spaces that there exists
a compactification $\overline{\clP}_N$ for $\clP_N$, classifying stable
curves with action of $G=\zmod{N}$. Let $k$ an algebraically closed field
of characteristic $p$ prime to $2N$ like in the previous subsection; here
we will just briefly find out the "cusps" of the geometric fibre
$\overline{\clP}_{N,[\chi]}=(\overline{\clP}_N\otimes k)_{[\chi]}$,
i.e. the points of the boundary.

We refer to~\cite{BR} for a precise definition of $\overline{\clP}_N$;
we will only need to know that there is a so-called "discriminant"
morphism
$$
\delta:\overline{\clP}_{N,[\chi]}\to \overline{\clM}_{0,(2,2)}
$$
with values in the stack of curves of genus 0 with four marked
points gathered by pairs. This morphism maps a Potts curve $C$
to the quotient $C/G$ marked by the branch points.

Recall the data $G=\zmod{N}$, $g=N-1$ and $\xi=\{\chi,\chi,\chi^{-1},\chi^{-1}\}$
of definition~\ref{def_Potts_curve}. In order to determine combinatorially
the stable curves lying on the boundary $\partial\overline{\clP}_{N,[\chi]}$
we use the combinatorial description of stable curves via their dual graph
$\Gamma$, as in~\cite{DM},~\cite{BR}. For a stable curve $C$, the graph
$\Gamma_C$ has the irreducible components of $C$ as vertices, and the double
points as edges.

\begin{prop}
The coarse moduli space of $\overline{\clP}_{N,[\chi]}$ is
$\overline{P}=\bbP^1\otimes k$, and the two cusps are topologically
the following two curves. The first has 2 branches isomorphic to
$\bbP^1$ intersecting $N$ times, and the second has 2 branches of
genus $\frac{N-1}{2}$ intersecting in only one point.

\begin{figure}[htb]
\centerline{\includegraphics{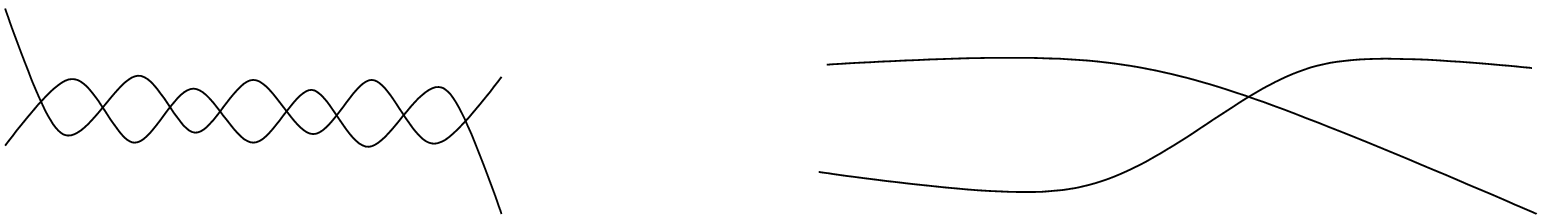}}
\end{figure}
\end{prop}

\begin{proo}
First, $\overline{P}$ is a normal proper curve of genus zero, hence it is
${\bbP}^1$. In dual graphs, we shall indicate marked points by wavy edges.
Hence let $(C,G)$
be a stable Potts curve. We have $\Gamma_\Sigma=\Gamma_C/G$, where $\Sigma=C/G$
has genus 0. We recall that $g'=\sum_i\,g_{\Sigma_i}+h^1(\Gamma_\Sigma)$, so that
the irreducible components $\Sigma_i$ of $\Sigma$ are all rational, and that
$\Gamma_\Sigma$ is a (connected) tree. Taking into account the four marked points
and the stability conditions, $\Gamma_\Sigma$ must be one of the two graphs:

\begin{figure}[htb]
\centerline{\includegraphics{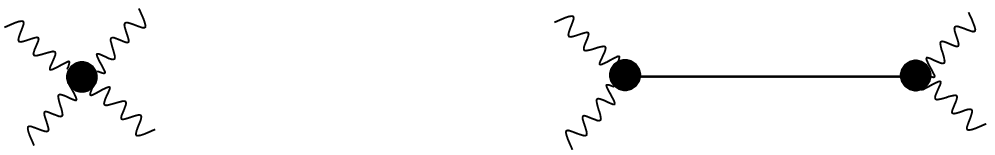}}
\end{figure}

But, $G$ being cyclic, a double point in $C$ maps to a double point
in $\Sigma$. Hence the first graph can not occur. So let $C_i$, $i=1,2$,
be irreducible components of $C$ above each of the components $\Sigma_i$
of $\Sigma$, and intersecting at a point $x$. The stabilizer of $x$
is $H=G_x$, and $h=[G:H]$ is its index. Take $a\in C_1$ a (smooth)
ramification point; the stabilizer of $C_1$ is $G_1=G$ because by
assumption $G=G_a\subset G_1$. Similarly
$G_2=G$, thus $C$ has only two irreducible components.

Now let us write the Riemann-Hurwitz formula for the quotient maps
$\pi|_{C_i}:C_i\to C_i/G_i\simeq\Sigma_i\simeq \bbP^1$. There is only
one orbit of double points, therefore their number is $[G:G_x]=h$:
$$
2g_i-2=N(-2)+\underbrace{2(N-1)}_{(r)}+\underbrace{h(N-h)}_{(d)}
=h(N-h)-2,\quad \forall i=1,2.
$$
where $(r)$ is the contribution of the ramification points, and $(d)$ the
contribution of the double points. In particular $g_1=g_2$, and we also know
that $g(C)=N-1=g_1+g_2+h^1(\Gamma_C)=g_1+g_2+h-1$. So, $h(N-h)=2g_1=N-h$, whence
$h=N$ or $1$.

Letting the family $y^2=x^{2N}+2x^N+t$ degenerate when $t\to 1$, we
get $y^2=(x^N+1)^2$. In this way we see that the first cusp described
above is the Potts curve of invariant $j=\infty$. Moreover, it is obvious
that there is only one (isomorphism class of) Potts curve with this
combinatorial aspect. A similar description with equations for the other
cusp would be more tricky. However, since we know that the other cusp can
not have the same combinatorics, we necessarily get the second picture
for $j=0$.
\end{proo}

\section{Moduli space of $\clP_p$}
\label{construction_EM_anychar}

\begin{theo*} \label{espace_modulaire_anychar}
The moduli space of $\clP_p$ is
${\bbA}^1_*\otimes\bbZ[\frac{\zeta+\zeta^{-1}}{2},\frac{1}{2}]$.
\end{theo*}

The rest of this section is devoted to the proof of the theorem. We keep
the former scheme of proof in two steps, but the complication coming from
wild ramification implies that we won't be able to "normalize" the
construction as well as before. Because of this we will need two additional
lemmas. In the first, which we now state, it is only for later convenience
that a primitive root of unity is denoted by $t$ instead of $\zeta$.

\begin{lemm*} \label{lemmPrelimInv}
Let $t,\psi$ be elements of a ring $A$, note $t^{[0]}=0$ and $t^{[i]}=1+t+\dots+t^{i-1}$
for $i\ge 1$. Let $\sigma$ be the endomorphism of the graded polynomial $A$-algebra
$A[X,Z]$ given by $\sigma(X)=t X+\psi Z$ and $\sigma(Z)=Z$.
Then $\sigma$ has exact order $p$ after any base change $A\to A'$ if and only
if $1+t+\dots+t^{p-1}=0$ and $(t-1,\psi)=A$. If this is so
then the algebra of invariants $A[X,Z]^\sigma$ is generated by $Z$ and the norm
$$
N(X,Z)=\prod_{i=0}^{p-1} \sigma^i(X)=\prod_{i=0}^{p-1} (X-t^{[i]}\psi Z)
$$
\end{lemm*}

\begin{proo}
It is clear that $Z$ plays no role, so we can dehomogenize and make $Z=1$.
Clearly, $\sigma^p=\Id \Leftrightarrow
t^p=1$ and $(1+t+\dots+t^{p-1})\psi=0$. Now let $A'$ run
through all residue fields $\kappa$ of $A$. Whenever $\overline{t}=1$ in
$\kappa$, the fact that $\sigma$ has exact order $p$ on the corresponding fibre
implies $\overline{\psi}\ne 0$. This means that $(t-1,\psi)=A$,
hence $1=u(t-1)+v\psi$ for some $u,v$.
From this we deduce that $1+t+\dots+t^{p-1}=0$. Conversely this
is easily seen to imply that $\sigma$ has exact order $p$ "universally".

As for the invariants we always have $A[N(X)]\subset A[X]^\sigma$. It is
that we have equality when $A$ is a field. In the general case, let $M$
be the cokernel of the inclusion, it is a finite $A$-module. The crucial
point is that, as the action is faithful fibrewise, the formation of
$A[X]^\sigma$ commutes with base change (this is a special case of results
concerning actions on smooth curves, see~\cite{BM2}, prop. 3.7). So for
every residue field we have $M\otimes \kappa=0$. By Nakayama's lemma,
it follows that $M=0$.
\end{proo}

\begin{rema*} \label{remaConducteurs_0}
We recall that for a $p$-Potts curve $(C,\sigma,\tau)$ over an algebraically closed field,
there are 2 fixed points in $C$ for the action of $\sigma$, and each has conductor
$m=1$. This follows from the description made in~\ref{CS}.
\end{rema*}

\titre{1st step: we build the morphism to the moduli space.}

Here we keep the notations of the proof of~\ref{espace_modulaire}. Let
$(f:C\to S,\sigma,\tau)$ be a $p$-Potts curve over $S$ in characteristic $p$. Let
$C\stackrel{r}{\to} E\stackrel{q}{\to} S$ be the factorization of $f$ through the
quotient by $\tau$. Then $\sigma$ induces an automorphism of order $p$ of $E$,
still denoted $\sigma$. By the remarks above, the divisor of fixed points $T=E^\sigma$,
finite of degree 2, is no longer \'etale but nevertheless fppf over $S$.
In fact, locally for the fppf topology it is a sum of two disjoint sections
$\Delta+\Delta'$, and above $S\otimes\bbF_p$ these sections have
the same support but infinitesimally they might be distinct (see below).

From now on we work with the section $\Delta$. Choosing ${\cal O}(1):={\cal O}(\Delta)$,
we set $V:=q_*{\cal O}(1)$ so $E=\bbP(V)$. The section corresponding to $\Delta$
gives a surjective map $h:V\to M=\clO(1)_{|\Delta}$, and unfortunately here we can not
go further to split the bundle. But $\ker(h)$ is known to be $N_\Delta^{_\vee}\otimes M$,
with $N_\Delta^{_\vee}$ the conormal sheaf of $\Delta$ in $E$: to see this,
remember the fundamental exact sequence
$$
0\to \Omega^1_{\bbP(V)/S}(1)\to q^*V\to \clO(1)\to 0
$$
and restrict to $\Delta$. Now
$N_\Delta^{_\vee}\simeq \clO(-\Delta)_{|\Delta}\simeq \clO(-1)_{|\Delta}$, so
that $\ker(h)\simeq \clO_S$. Hence we have an extension
$$
0\to \clO_S\to V\to M\to 0
$$
Let us see now how $\sigma$ acts on this. As an automorphism of $E$,
it pulls back ${\cal O}(1)$ to an invertible sheaf of degree 1, i.e. there is
an isomorphism $u_\sigma:\sigma^{*}\clO(1) \simeq q^{*}K\otimes \clO(1)$
for some $K\in\Pic(S)$. Moreover $\sigma$ is the identity on $\Delta$, so that
restricting $u_\sigma$ to $\Delta$ shows that $K$ is trivial. 
Now $\sigma$ is given by a surjective morphism of sheaves
$$q^*V\,\to\,\sigma^*{\cal O}(1)\stackrel{u_\sigma}{\simeq}{\cal O}(1)$$
Taking direct images by $q$, we obtain an automorphism
$\varphi_\sigma:V\to V$. This map induces an automorphism of
$M$ and of $\ker(h)\simeq \clO_S$, and is well determined up to an invertible
global section of $\clO_S$, but requiring ${\varphi_\sigma}|_{\clO_S}={\rm id}$
makes $\varphi_\sigma$ canonical. Now the action on $M$ is given by
multiplication by a global section $t\in \Gamma(S,\clO_S^\times)$,  this
means that if we choose locally a coordinate $X$ for $M$ as in the proof
of~\ref{espace_modulaire}, the action is
$$
\begin{array}{ll}
\varphi_\sigma(X)=t X+\psi Z & (\mbox{for some } \psi\in \Gamma(M,\clO_S)) \\
\varphi_\sigma(Z)=Z &
\end{array}
$$
As $\sigma$ has order $p$ on the fibres, we have $1+t+\dots+t^{p-1}=0$,
and $t-1$ and $\psi$ generate $\clO_S$, by lemma~\ref{lemmPrelimInv}.
In particular, denoting by $\omega= t^{[1]}\dots t^{[p-1]}$
we have $p=\omega (t-1)^{p-1}$ (notations of lemma~\ref{lemmPrelimInv}).

The double cover $r$ is described by the decomposition
$r_*\clO_C=\clO_E\oplus \clL$ and by a section $\theta\in\Gamma(E,\clL^{ -2})$,
well determined up to an element of $\Gamma(E,\clO_E^\times)$. Also,
$\clL^{ -2}\simeq {\cal O}_E(2p)$, so we can identify $\theta$ with a
$\sigma$-invariant section of
$\Gamma(E,{\cal O}_E(2p))=\oplus_{j=0}^{2N}\,\Gamma(S,M^{ j})$.
By lemma~\ref{lemmPrelimInv} we get
$$\theta\propto H=U N(X,Z)^2+A N(X,Z) Z^p+BZ^{2p}$$
for some sections $U,A,B$ of $\clO_S,M^p,M^{2p}$. As in the proof
of~\ref{espace_modulaire} we see that $U$ is invertible. Now we must express
that the fibres of $C/S$ are smooth, i.e. that $\theta$ has no multiple zero.
To this aim we compute its discriminant, it turns out that
$$
\Resul(H,H')=-U^p\,\omega^{2p}\,\delta^{p-1}\,
(A^2-4BU)^p
$$
where $\delta:=H(-\psi,t-1)=U\psi^{2p}-A\psi^p(t-1)^p+B(t-1)^{2p}$ and we
recall that $\omega= t^{[1]}\dots t^{[p-1]}$. In this expression both $U$
and $\omega$ are invertible. So the smoothness condition is that $\Resul(H,H')$,
or equivalentlty $\delta(A^2-4UB)$, is invertible (remark : this contains
the condition that $t-1$ and $\psi$ generate $\clO_S$). We are led to define
$$j=\frac{U(U\psi^{2p}-A\psi^p(t-1)^p+B(t-1)^{2p})}{A^2-4UB}$$
It lies in $\Gamma(S,\clO_S^\times)$ because numerator and denominator
are invertible sections of $M^{2p}$. Also it is independant of the leading
coefficient of $H$. If our invariant were merely the discriminant of $H$
then, being intrisic, it would be obviously invariant under changes of variables.
Here this is not the case, so we will proceed to check this in the form of
a lemma.

\begin{lemm*}
The definition of $j$ is independant of the choice of coordinate $X,Z$
and of the choice of section $\Delta$.
\end{lemm*}

\begin{proo}
From now on we always normalize the expressions by setting $U=1$.
First, in the choice of the coordinate system the only loose variable
is $X$ so we can assume that $Z'=Z$ and $X'=\alpha X+\beta Z$. Then
$\sigma(X')=t X'+\psi'$ where $\psi'=\alpha\psi-(t-1)\beta$. The norm
$N'(X',Z')$ with $t'=t$ and $\psi'$ is
$$
N'(X',Z')=\prod_{i=0}^{p-1} (X'-t^{[i]}\psi' Z')
=\prod_{i=0}^{p-1} (\alpha X+\beta Z-t^{[i]}\psi' Z)
$$
Being $\varphi_\sigma$-invariant, this is a polynomial in $N$ and $Z^p$, hence
there exists $\xi$ such that $N'(X',Z')=\alpha^pN(X,Z)+\xi Z^p$. The polynomial
$H'=N'(X',Z')^2+A' N'(X',Z') Z'^p+B'Z'^{2p}$ associated to $\theta$ can be
expressed in terms of $N(X,Z)$ and $Z^p$ (here $H'$ is not the derivative of $H$~!).
As the change of variables $(N,Z^p)\leftrightarrow (N',Z'^p)$ has determinant
$\alpha^p$, the discriminant of $H'$ viewed as a polynomial of degree~2 in
$(N,Z^p)$ is $\alpha^{2p}(A'^2-4B')$.
Of course since $H\propto H'$ we have $H'=\alpha^{2p}H$. Computing discriminants
gives $\alpha^{2p}(A'^2-4B')=\alpha^{4p}(A^2-4B)$. Then substituting
$(X,Z)=(-\psi,t-1)$, first in $N'=\alpha^pN+\xi Z^p$, second in $H'=\alpha^{2p}H$,
gives $\delta'=\alpha^{2p}\delta$. Finally $j'=\delta'/(A'^2-4B')=\delta/(A^2-4B)=j$.

Now, assume that we choose the section $\Delta'$ instead of
$\Delta$. It is not as obvious as in~\ref{espace_modulaire} that formally
this has the effect of changing $t$ to $t^{-1}$, and not even that the invariant
will not change; so we sketch the details. The equation of $\Delta'$
is given by the new coordinate $Z'=\varphi_\sigma(X)-X=(t-1)X+\psi Z$.
As we saw just above, we can choose $X'$ as we like; to simplify matters
we recall that there exist sections $u,v$ such that $u(t-1)+v\psi=1$ and
we choose $X'=vX-uZ$ so as to have a unimodular change of variables
$$
\left\{
\begin{array}{l}
X'=vX-uZ \\
Z'=(t-1)X+\psi Z \\
\end{array}
\right\}
\qquad\Leftrightarrow\qquad
\left\{
\begin{array}{l}
X=\psi X'+uZ' \\
Z=-(t-1)X'+vZ' \\
\end{array}
\right\}
$$
Then we have $\varphi_\sigma(X')=X'+vZ'$ and $\varphi_\sigma(Z')=t Z'$.
The condition that $\varphi_\sigma$ acts trivially on $Z'$ leads to consider
$\varphi_\sigma'=t^{-1}\varphi_\sigma$, and we obtain
$$
\begin{array}{l}
\varphi_\sigma'(X')=t^{-1}(X'+vZ') \\
\varphi_\sigma'(Z')=Z'
\end{array}
$$
With $t'=t^{-1}$ and $\psi'=t^{-1}v$, the norm is
$N'(X',Z')=\prod_{i=0}^{p-1} (X'-(t^{-1})^{[i]}t^{-1}v Z')$. Using
that $-(t^{-1})^{[i]}t^{-1}=t^{[p-i]}$ we compute
$$
N'(X',Z')=\prod_{i=0}^{p-1} (vX-uZ-t^{[i]}Z)
$$
which is $\varphi_\sigma$-invariant, so there exists $\xi$ such that
$N'(X',Z')=v^pN(X,Z)+\xi Z^p$. Substituting $(X,Z)=(-\psi,t-1)$
yields immediately $-1=-v^p\psi^p+(t-1)^p\xi$, so the change of
variables between the invariants of degree $p$ is unimodular :
$$
\left\{
\begin{array}{l}
Z'^p=\psi^p Z^p+(t-1)^p N(X,Z) \\
N'(X',Z')=\xi Z^p+v^pN(X,Z) \\
\end{array}
\right.
$$
Thus for $H'=N'(X',Z')^2+A' N'(X',Z') Z'^p+B'Z'^{2p}$, its discriminant
as a polynomial of degree~2 in $(N,Z^p)$ is still $A'^2-4B'$. Now, when
expressed in terms of $(N,Z^p)$ the polynomial $H'$ has leading coefficient
$\delta'=H'(-t^{-1}v,t^{-1}-1)$, so $H'=\delta'H$. Computing discriminants
gives $A'^2-4B'=\delta'^2(A^2-4B)$, and substituting $(X,Z)=(-\psi,t-1)$
gives $1=\delta'\delta$, so here again $j'=j$.
\end{proo}

Therefore the only change is that having chosen $\Delta'$ instead of
$\Delta$, we recover $t^{-1}$ instead of $t$ as a $p$-th root of unity.
So we have a well-defined map
$$
\bbZ[\frac{\zeta+\zeta^{-1}}{2},\frac{1}{2}]
\left[X,X^{-1}\right]\to \Gamma(S,{\cal O}_S)^\times
$$
if we set $F(\zeta+\zeta^{-1})=t+t^{-1}$ and $F(X)=j$.
Eventually we have a map
$S\to {\bbA}^1_*\otimes\bbZ[\frac{\zeta+\zeta^{-1}}{2},\frac{1}{2}]=P$.
In fact we completed the construction after base change to $T=E^\sigma$,
but the construction being canonical, by fppf descent we obtain a morphism
defined on $S$. It is clear that the construction is functorial
in $S$, providing a morphism $\Phi:\clP_p\to P$.

\titre{2nd step: we check the properties of a moduli space.}

Here, provided we give a family of Potts curves with a finite,
surjective associated morphism to the moduli space, the arguments
of the proof of~\ref{espace_modulaire} carry on. We consider the norm
for $\psi=1$, $N(x)=\prod_{i=0}^{p-1} (x-\zeta^{[i]})$ and the curve
with equation $y^2=N(x)^2+\lambda N(x)+1$, with invariant
$j_0(\lambda)=\frac{1-\lambda(\zeta-1)^p+(\zeta-1)^{2p}}{\lambda^2-4}$,
over the base
$$
S_0=\Spec\left(\bbZ[\zeta,\frac{1}{2}]
\left[\lambda,j_0(\lambda),\frac{1}{j_0(\lambda)}\right]\right)
$$
As in the proof of~\ref{espace_modulaire}, to be rigorous we can easily
give another smooth affine part for this curve, but we omit this detail.
Then the proof of~\ref{espace_modulaire} works similarly, ending the
proof of th.~\ref{espace_modulaire_anychar}. \qed

\section{The fibre of $\clP_p$ at the prime $p$} \label{MSWC}

At last we study the stack of $p$-Potts curves in characteristic $p$,
that is to say $\clP_p\otimes \bbF_p$. First we compute its moduli space,
along the same lines as above; the main difficulty here is that this space
is not normal anymore, so we need the help of deformation theory as well as
the operation of "2-quotient" of an algebraic stack (see~\cite{Ro}, \cite{AbV}).
The result shows that the moduli space of the $\bbZ[1/2]$-stack $\clP_p$ has good
reduction at $p$. Then, we compute the Picard group of $\clP_p\otimes \bbF_p$.

\subsection{Moduli space in characteristic $p$}
\label{EMCPp}

We still assume that $p>3$. Before we state the theorem, let us consider the
problem of constructing the moduli space of the stack $\clP_p\otimes \bbF_p$
along the same lines as above (section~\ref{construction_EM_anychar}).
In the first step of the proof of~\ref{espace_modulaire_anychar} nothing needs
any change (except perhaps the fact that $\psi$ can be chosen to be equal to 1;
this makes the computations slightly simpler but is anecdotical). We obtain
a classifying morphism $\Phi$
from $\clP_p$ to $P=\bbA^1_*\otimes\bbF_p[\frac{\zeta+\zeta^{-1}}{2}]$.
Here, throughout the construction, $\frac{\zeta+\zeta^{-1}}{2}$ is still a root
of the polynomial $\psi_p$ of~\ref{remaMatrix}, which is none other than
the polynomial $\psi$ of~\cite{BM1}, 4.2.5 and 4.2.6. In characteristic $p$
we simply have $\psi(X)=X^{(p-1)/2}$, so
$P={\bbA}^1_*\otimes \frac{\bbF_p[z]}{z^{(p-1)/2}}$.

The difference comes with the 2nd step where we must check the properties
of a moduli space. The problem is that $P$ is not a normal scheme anymore,
so we will have to use a different strategy. One ingredient will be the
"2-quotient" of an algebraic stack whose objects all possess a fixed finite
group inside their automorphism group. Here is a brief summary of its
properties (\cite{Ro}, chap. I, prop. 3.0.2 or \cite{AbV}, prop. 3.5.1):

\begin{prop} \label{Twoquotient}
Let $S$ be a scheme and $\clM$ an algebraic stack over $S$. Let $G$ be a finite
group, assume that for every object $x\in \clM(T)$ there is an injection
$i_x:G_T\hookrightarrow \Aut_T(x)$, whose formation is compatible
with cartesian diagrams (see~{\rm \cite{Ro}}). Then there exists an algebraic stack
$\clM\twoquotient G$ and a map $f:\clM\to \clM\twoquotient G$, such that $f$ maps the
elements of $G$ to the identity, and is universal with respect to this property.
The geometric points of $\clM\twoquotient G$ are the same as those of $\clM$,
but for $x$ such a point we have $\Aut_{\clM\twoquotient G}(x)=\Aut_\clM(x)/G$.
The formation of $\clM\twoquotient G$ commutes with base change on $S$. Moreover,
$f$ is an \'etale gerbe. The stack $\clM\twoquotient G$ has a coarse moduli space
if and only if $\clM$ has one, and if this is the case the moduli spaces are
the same. Finally, if $\clM$ is separated or proper, then $\clM\twoquotient G$
has the same properties. \qed
\end{prop}

A basic example of this is given by the classifying stack $BG$ of a finite
abelian group $G$ over a scheme $S$. Its objects are $G$-torsors, and $G$ lies
in all automorphism groups. In this case $\clM=BG$, and $\clM\twoquotient G=S$.
Now let us come back to the stack $\clP_p\otimes \bbF_p$. We are going to show~:

\begin{theo} \label{espace_modulaire_carp}
If $G=\zmod{2}\times\zmod{p}$, then we have an isomorphism
$$
(\clP_p\otimes \bbF_p)\twoquotient G \isomto {\bbA}^1_*\otimes \frac{\bbF_p[z]}{z^{(p-1)/2}}
$$
In particular, ${\bbA}^1_*\otimes \frac{\bbF_p[z]}{z^{(p-1)/2}}$
is the coarse moduli space of $\clP_p\otimes \bbF_p$.
\end{theo}

From proposition~\ref{automorphismes_car_p} we know that $p$-Potts curves
(over any base $S/k$) have the constant group scheme $G=\zmod{2}\times\zmod{p}$
as their automorphism group (use non-ramification of $\Aut_S(C,\sigma,\tau)$).
Hence by definition of the 2-quotient the morphism $\Phi$ factors through
$$
\Psi:(\clP_p\otimes \bbF_p)\twoquotient G\to
{\bbA}^1_*\otimes \frac{\bbF_p[z]}{z^{(p-1)/2}}
$$
with the stack $\clQ:=\clP_p\otimes \bbF_p\twoquotient G$ representable
by an algebraic space (\cite{LMB}, 8.1.1). Thanks to deformation theory,
known from~\cite{BM1}, we shall show that $\Psi$ is \'etale~:

\begin{prop} \label{eq_def_univ2}
Let $(C,\sigma,\tau)$ be a $p$-Potts curve over $k$.
Then the ring that prorepresents the functor of deformations of $C$ to local,
artinian $k$-algebras is $\frac{\bbF_p[z]}{z^{(p-1)/2}}[\![t]\!]$.
\end{prop}

\begin{proo}
This is an example of the computation of deformations of $\zmod{p}$-actions
on smooth curves by Bertin and M\'ezard~\cite{BM1}. In their work everything
is done over an algebraically closed field, as is usual in arithmetic
geometry in mixed characteristic, but one can check that this assumption
is not necessary in their article. Indeed the results they use are Schlessinger's
criteria that need no assumption on the base field, and Serre's \emph{Corps
Locaux} that uses perfect fields. So their results are valid over $\bbF_p$.

Now, as $\tau$ has order 2 which is assumed to be prime to $p$, the deformation
ring of $(C,\sigma,\tau)$ is the deformation ring of $(D=C/\tau,\sigma)$ (we
still denote $\sigma$ the automorphism induced on $C/\tau$). We first look
at deformations of $(D,\sigma)$ to algebras over the ring of Witt vectors
$W(\bbF_p)$, like in~\cite{BM1}. In the article, corollary~3.3.5 shows that 
the universal deformation ring of $(D,\sigma)$ is
$$
R_{\rm gl}=(R_1\widehat{\otimes}\dots\widehat{\otimes} R_r)[\![U_1,\dots,U_N]\!]
$$
with $N={\rm dim}_k\,H^1(D/\sigma,\pi_*^\sigma({\cal T}_D))$. Here
there is only $r=1$ orbit of fixed points, with conductor $m=1$ (see
remark~\ref{remaConducteurs_0}). Also $R_1=W(\bbF_p)[\![X]\!]/\psi(X)$
by~\cite{BM1}, theorem 4.2.8. Finally, the computation of $N$ is done in the
course of the proof of theorem 4.2.8, namely $N=1$. Reducing modulo $p$, we
have $R_1/pR_1=\frac{\bbF_p[z]}{z^{(p-1)/2}}$. We obtain the universal ring for
deformations to $k$-algebras as $\frac{\bbF_p[z]}{z^{(p-1)/2}}[\![u]\!]$.
\end{proo}

\noindent {\bf End of proof of theorem~\ref{espace_modulaire_carp}~:}
As the 2-quotient $\clP_p\otimes \bbF_p\to \clQ$ is \'etale, it follows from
the proposition that the extensions of complete local rings corresponding to
$\Psi$ are trivial, meaning that both $\Phi$ and $\Psi$ are \'etale. A first
consequence is that $\clQ$ is not only an algebraic space, but in fact a
scheme. Also, considering the scheme $\clQ'$ defined by the fibre square:
\begin{diagram}[height=2em,width=2em,PostScript=dvips]
\clQ' & \rTo^{u'} & {\bbA}^1_*\otimes k \\
\dTo  &           & \dTo \\
\clQ  & \rTo^u    & {\bbA}^1_*\otimes \frac{\bbF_p[z]}{z^{(p-1)/2}}
\end{diagram}
we have that $u'$ is \'etale, hence $\clQ'$ reduced. Moreover $u'$ is bijective,
so by Zariski's Main Theorem for schemes $u'$ is an isomorphism. So $\clQ_\red=\clQ'$
is affine, which implies that $\clQ$ itself is. Now using \cite{SGA1}, Exp. I, th. 6.1
we get that $u$ is an isomorphism. \qed

\subsection{The Picard group of $\clP_p\otimes \bbF_p$} \label{GPPp}

Let $k$ be an algebraically closed field of characteristic $p\ne 2$ (actually the assumption
of algebraic closure will not be necessary). Let $P$ be the moduli space of
$\clP_p\otimes k$ and $A=\frac{k[z]}{z^{(p-1)/2}}[X,\frac{1}{X}]$ its ring of functions.
The adaptation of the arguments of~\ref{GPP} gives the following result~:

\begin{theo} \label{PicardWild}
There is an isomorphism $\Pic(\clP_p\otimes k)\simeq \zmod{2}\times (1+zA)$
where $1+zA$ is a subgroup of the multiplicative group of invertible elements
in $A$.
\end{theo}

\begin{proo}
Put $\clP:=\clP_p\otimes k$. Let $U\to \clP$ an open set of the \'etale
site of $\clP_p$, corresponding to a $p$-Potts curve $(C,\sigma,\tau)$
over $U$. We know that $\Aut_U(C,\sigma,\tau)$ is the constant group scheme
$G=\zmod{2}\times\zmod{p}$. Then, as in~\ref{LemmBeta}, for any invertible
sheaf $L$ and any $g\in G$ we have an automorphism $\ell_{\Id,g}:L|_U\isomto L|_U$
given by a global section of $\clO_U^\times$. This gives a map $\beta$ from
$\Pic(\clP)$ to the abelian group of homomorphisms
of abelian sheaves over $\clP$ from $G$ to
$\clO_{\clP}^\times$. Let $q:\clP\to P$ be
the map to the moduli space, then we have an exact sequence:
$$
0\to \Pic(P)\stackrel{q^*}{\longrightarrow}
\Pic(\clP) \stackrel{\beta}{\longrightarrow}
\Hom_{\clO_{\clP}}(G,\clO_{\clP}^\times)\to 0
$$
Indeed, exactness in the middle is proved exactly by the proof
of~\ref{PicardTame}. The pullback $q^*$ is injective because
$\Pic(P)=0$. The fact that $\beta$ is surjective is also clear because,
given a character $f:G\to \clO_{\clP}^\times$,
we can twist the structure sheaf $\clO_{\clP}$
so as to define a sheaf $L$ by $L|_U=\clO_U$ for all $U$, and
$\ell_{\Id,g}:\clO_U\isomto \clO_U$ equal to multiplication by
$f(U)(g)$. This sheaf satisfies $\beta(L)=f$.

It remains to compute $\Hom(G,\clO_{\clP}^\times)$, which is just
the character group of $G$. Using adjunction and the property of the
moduli space that $q_*\clO_{\clP}=\clO_P$, we have
$$
\Hom_{\clP}(G,\clO_{\clP}^\times)
=\Hom_P(G,\clO_P^\times)=\mu_2(A)\times\mu_p(A)=\zmod{2}\times (1+zA)
$$
and this is the result.
\end{proo}

%

\end{document}